\documentclass{aims_noheader}
\usepackage{amsmath}
\usepackage{amssymb}
  \usepackage{paralist}
  \usepackage{graphics} 
  \usepackage{epsfig} 
\usepackage{graphicx}  \usepackage{epstopdf}
 \usepackage[colorlinks=true]{hyperref}
\usepackage{todonotes}
\hypersetup{urlcolor=blue, citecolor=red}

\def\currentvolume{}
 \def\currentissue{}
  \def\currentyear{}
   \def\currentmonth{}
    
\newtheorem{theorem}{Theorem}[section]

\theoremstyle{definition}

\newtheorem{remark}{Remark}

\title[ABC for the Westervelt equation] 
       {Absorbing boundary conditions for the Westervelt equation}
      
\author[Barbara Kaltenbacher and Igor Shevchenko]{}
                            


\email{i.shevchenko@imperial.ac.uk}
\email{Barbara.Kaltenbacher@uni-klu.ac.at}
 
\thanks{The first author is supported in part by EPSRC Mathematics Platform grant EP/I019111/1}

\begin{document}
\maketitle

\centerline{\scshape Barbara Kaltenbacher}
\medskip
{\footnotesize
 \centerline{Alpen-Adria-Universit\"at Klagenfurt}
   \centerline{Institute of Mathematics}
   \centerline{Klagenfurt, A-9020, Austria}
}

\medskip

\centerline{\scshape Igor Shevchenko}
\medskip
{\footnotesize
 \centerline{Imperial College London}
   \centerline{Department of Mathematics}
   \centerline{London, SW7 2AZ, UK}
} 

\bigskip


\begin{abstract}
The focus of this work is on the construction of a family of nonlinear absorbing boundary
conditions for the Westervelt equation in one and two 
space dimensions. The principal ingredient used in the design of such conditions is pseudo-differential calculus. This
approach enables to develop high order boundary conditions in a consistent way which are 
typically more accurate than their low order analogs. Under the hypothesis of 
small initial data, we establish local well-posedness for the Westervelt equation with the absorbing boundary conditions. 
The performed numerical experiments illustrate the efficiency of the proposed boundary
conditions for different regimes of wave propagation.
\end{abstract}

\section{Introduction}
Constantly growing needs of numerical simulations in science and engineering
often require considering problems which are naturally formulated in unbounded domains.
Typical examples can be found in many problems originating from fluid dynamics, solid
mechanics, aerodynamics, electrodynamics, acoustics, etc.
However, the numerical solution of such problems requires a finite region. There are basically 
two approaches which can be used to reformulate problems in infinite domains as
problems in finite domains. The first approach is to map the originally unbounded domain to 
a bounded one. Simple as the problem sounds the solution in practical applications is far from known.
This is mostly due to reasons which are connected with singularities of
the new equation that results from the mapping.
The second approach, which we follow in this work, is to impose
fictitious boundaries to truncate the domain of interest. Such artificial boundaries require special boundary conditions 
so that the boundary value problem is well-posed and its solution is an accurate approximation to the restriction of the solution
in the unbounded domain. In other words, these boundary conditions have to be transparent to or, as they are usually called, absorbing for 
solutions propagating outwards the artificial boundary. 
It is commonly recognized that absorbing boundary conditions (ABCs)
play a key role in computations on unbounded domains and have a profound impact on the accuracy of numerical methods.
Over the past thirty years, ABCs have developed into a vigorous 
research direction including a wide spectrum of methods and approaches. The description of these techniques is out of the 
scope of this work and therefore we restrict ourselves to referring the reader to the comprehensive review
articles~\cite{Givoli1991,Tsynkov1998,Hagstrom1999,Hagstrom2003,Givoli2004,GivoliBookChapter08} and the references therein. 

The focus of this work is on construction of ABCs for high-intensity focused ultrasound (HIFU) which plays an important role
in many medical and industrial applications such as diagnostic ultrasound~\cite{FlochFink1997,SimonVanBarenEbbini1998,PernotWatersBercoffTanterFink2002},
thermotherapy of tumors~\cite{FlochTanterFink1999,HallajClevelandHynynen2001,ConnorHynynen2002},
lithotripsy~\cite{AverkiouCleveland1999}, ultrasound cleaning and sonochemistry. Linear
models of wave propagation are not applicable in HIFU due to nonlinear effects requiring more sophisticated acoustic equations to be taken into account.
In this work, we develop local in space and time ABCs for the Westervelt equation used as a basic 
acoustic model in various HIFU simulations. The Westervelt equation is one of the fundamental equations
governing the propagation of acoustic waves in nonlinear regimes~\cite{Westervelt1963,HallajClevelandHynynen2001,ConnorHynynen2002,Clason2009}:
\begin{equation}
u_{tt}-c^2\Delta u-b\Delta u_t=\frac{\beta_a}{\varrho c^2}(u^2)_{tt}\quad\text{in }(0,T)\times\Omega,
\label{Westervelt0}
\end{equation}
where $\Omega\subseteq\mathbb{R}^d$, $d\in\{1,2,3\}$, 
$u=u(\cdot,t)$ is the acoustic pressure, $c>0$ is the speed of sound, $b>0$ is the acoustic diffusivity, $\varrho>0$ is the mass density,
$\beta_a=1+B/(2A)$ with $B/(2A)>0$ standing for the parameter of nonlinearity of the fluid, 
$T$ is the final time at which the problem is to be solved. All the parameters are assumed to be constant.
We rewrite \eqref{Westervelt0} in a form more convenient for further treatment
\begin{equation}
c^{-2}u_{tt}-\Delta u-\beta\Delta u_t=\gamma(u^2)_{tt}\quad\text{in }(0,T)\times\Omega
\label{Westervelt}
\end{equation}
with $\beta=b/c^2$, $\gamma=\beta_a/(\varrho c^4)$, and complement \eqref{Westervelt} by initial conditions 
\begin{equation}
u(\cdot,t=0)=u_0\,, \quad u_t(\cdot,t=0)=u_1\quad\text{in }\Omega,
\label{IC}
\end{equation}
and by inhomogeneous Neumann and absorbing boundary conditions
\begin{equation}
u_n\Big|_{(0,T)\times\Gamma_{\rm N}}=g(t), \quad \mathcal{A}u\Big|_{(0,T)\times\Gamma_{\rm A}}=0,
\label{BC}
\end{equation}
where $\partial\Omega=\Gamma_{\rm N}\cup\Gamma_{\rm A}$, subscript $n$ denotes the normal derivative on the boundary, 
and the operator $\mathcal{A}$, on the absorbing boundary $\Gamma_{\rm A}$, is an annihilating operator for outgoing waves which we specify in due course.

In spite of the intensive research activity in the field of transparent boundary conditions, most results have been obtained for linear problems with
constant coefficients. Wave equations with variable coefficients have received much less attention, not to mention nonlinear models.
There are only few papers devoted to problems with variable coefficients~\cite{EngquistMajda1979}, 
convective~\cite{BecacheGivoliHagstrom2010} and nonlinear~\cite{Hedstrom1979,Szeftel2006_2,ZhangXuWu2009,PazStortiGarelli2009} terms.
Despite the existence of some approaches to the construction of ABCs for nonlinear wave models their application to concrete equations
is rather sophisticated and still out of the scope of most research works.

In this work we design ABCs based on the theory of pseudo-differential~\cite{KohnNirenberg1965,Hormander1965,Nirenberg1975} calculus. We will also address a possible approach via para-differential~\cite{Bony1981,Meyer1981} calculus in the appendix.
The first approach is applicable to linear wave equations with variable coefficients. Therefore it is used for the Westervelt equation linearized in a neighborhood of a reference solution. 
The second approach will be directly applied to
the nonlinear Westervelt equation. Before going into detail with the derivation of ABCs, we remark that both theories have been already used in the construction of transparent boundary conditions.
For example, the pseudo-differential calculus was exploited by Engquist and Majda in~\cite{EngquistMajda1979} to design ABCs for the linear wave equation with variable coefficients.
Transparent boundary conditions for the semilinear wave equation as well as for the nonlinear Schr\"odinger equation with the help of para-differential operators were obtained 
in~\cite{Szeftel2006_1} and~\cite{Szeftel2006_2}, respectively.

\begin{remark}\label{rem:beta0}
The nonlinearity in the Westervelt equation comes along with a strong damping term $b\Delta u_t$. In fact, this strong damping, besides being a physically imposed term, also plays a quite particular mathematical role. As already observed in \cite{KL09}, \cite{KaltenbacherLasiecka2010} in the context of different boundary conditions, strong damping $\beta>0$ is essential in two and higher space dimensions in order to compensate the nonlinearity and avoid degeneracy in the equation, (whereas in 1-d we have wellposedness also in case $\beta=0$).
On the other hand, the strong damping term destroys the wave like character of the equation since it implies decay of the energy and a rather parabolic than hyperbolic behaviour of the equation, cf. \cite{KL09}. This results in the observation that the (linearized) differential operator defining the Westervelt equation with strong damping is not amenable to a factorization (see, e.g., \eqref{NFact_Westervelt_1D_u1} below) as required for constructing absorbing boundary conditions. For this reason, we skip the strong damping term during derivation of the ABCs. Of course it has to stay in the PDE, though (for physical reasons and since otherwise wellposedness would fail, as mentioned above). It is clear that the inevitable use of integration by parts in deriving energy estimates for the PDE causes the appearance of boundary terms resulting from the presence of the term $b\Delta u_t$. So these terms finally have to be taken into account in the ABC as well. As mentioned already, the factorization approach based on pseudo- or paradifferential calculus is not appropriate for doing so. Thus, incorporation of the $\beta$-term in the boundary conditions will be done as a postprocessing step after the pseudodifferential factorization, and it will be done on the basis of energy considerations. The latter will also allow to prove well-posedness of the resulting initial boundary value problems for the Westervelt equation.
\end{remark}

\subsection{Main results} 

The novelty of our work lies in the derivation and analysis of high-order ABCs for the Westervelt equation which have not been construct so far. 
We will do so for the one- and two dimensional versions of the Westervelt equation~(\ref{Westervelt}) first of all in a domain without corners, see Section \ref{sec:ABCWestervelt}. Additionally, we will 
provide well-posedness results for the Westervelt equation with zero and first order conditions in one and two space dimensions in Section \ref{secwell}.

In this section we summarize the boundary conditions derived in this paper together with the main well-posedness results.

\medskip

For the case of one space dimension we will derive zero, 
\begin{equation}
\left.\left(u_n+\sqrt{c^{-2}-2\gamma u}\ u_t\right)\right|_{\partial\Omega} =0,
\label{ABC0_1D_intro}
\end{equation}
first 
\begin{equation}
\left.\left(u_n+\sqrt{c^{-2}-2\gamma u}\ u_t-\frac{\gamma}{2\sqrt{c^{-2}-2\gamma u}}\left(
u_t u -\frac{1}{\sqrt{c^{-2}-2\gamma u}} u_n u\right)\right)
\right|_{\partial\Omega}=0 
\label{ABC1_1D_intro}
\end{equation}
and second order ABC 
\begin{equation}
\left.\left(u_{nt}+\sqrt{c^{-2}-2\gamma u}\, u_{tt}-\frac{\gamma}{2\sqrt{c^{-2}-2\gamma u}}\left(
(u_t)^2 -\frac{1}{\sqrt{c^{-2}-2\gamma u}}\, u_n u_t-\mu(u) u\right)\right)
\right|_{\partial\Omega}=0 
\label{ABC2_1D_intro}
\end{equation}
(with $\mu$ defined as in \eqref{mu}) in Section \ref{sec:LinWestervelt_1d}.

In 2-d, the zero order ABCs derived in Section \ref{sec:LinWestervelt_2d} are 
\begin{equation}
\left.\left(u_n+\sqrt{c^{-2}-2\gamma u} u_t\right)\right|_{\partial\Omega} =0
\label{ABC0_2D_intro}
\end{equation}
and the first order ones are 
\begin{equation}
\begin{aligned}  
&\left((u_{nt}+\sqrt{c^{-2}-2\gamma u}\, u_{tt}-\frac{1}{2\sqrt{c^{-2}-2\gamma u}} u_{\vartheta\vartheta}
\right.\\
&-\frac{\gamma}{2\sqrt{c^{-2}-2\gamma u}}
\left(u_t -\frac{1}{\sqrt{c^{-2}-2\gamma u}} u_n \right)u_{t}\\
&\left.\left.-\frac{\gamma}{2(c^{-2}-2\gamma u)^{3/2}}
\left(\frac12 u_t -\frac{1}{\sqrt{c^{-2}-2\gamma u}} u_n \right)\int_0^\cdot u_{\vartheta\vartheta}\, dt\right)
\right|_{\partial\Omega}=0\,,
\end{aligned} 
\label{ABC1_2D_intro}
\end{equation}
where subscript $\vartheta$ denotes the tangential derivative.
Here the zero order ABC and the first line of the first order ABC are exactly what one would expect from the linear case with constant coefficients. 

Note that for reasons outlined in Remark \ref{rem:beta0} above we set $\beta=0$ in these derivations.
The energy considerations in Section \ref{subsec:beta} will allow us to appropriately take into account the third order derivative term going with $\beta$. With the according modifications in the ABCs \eqref{ABC0_1D_intro}, \eqref{ABC1_1D_intro}, \eqref{ABC0_2D_intro}, \eqref{ABC1_2D_intro}, and denoting 
$$ u^\beta:= u+\beta u_t$$
we will obtain the following local in time well-posedness results for sufficiently small initial data $u_0,u_1$, and 
\[
u_2=\frac{1}{c^{-2}-2\gamma u_0}\left(\Delta u_0+\beta \Delta u_1+2\gamma (u_1)^2\right)
\]

\begin{theorem}\label{theo:zero1d}
For $\beta\ge0$, any open interval $\Omega=(a,b)\subseteq\mathbb{R}$ and any $T>0$ there exists $\rho>0$ such that for all 
initial data $u_0,u_1$ satisfying 
$\|u_2\|_{L^2(\Omega)}+\|u_1\|_{H^1(\Omega)}<\rho$,
a solution $u\in C^2(0,T;L^2(\Omega))\cap C^1(0,T;H^1(\Omega))$ to 
\begin{equation}
\begin{split}
& (c^{-2}-2\gamma u)u_{tt}-u_{xx}-\beta u_{txx}=2\gamma (u_t)^2\quad\text{in } (0,T)\times\Omega,\\
& u(t=0)=u_0\,,\ u_t(t=0)=u_1\quad\text{in }\Omega,\\
& u^\beta_n+\sqrt{c^{-2}-2\gamma u}\, u_t=0\quad\text{at } (0,T)\times\{a,b\}
\end{split}
\label{Westervelt_zero_1D_IBVP_intro}
\end{equation}
exists and is unique.
\end{theorem}

\begin{theorem}\label{theo:first1d}
For $\beta\ge0$, any open interval $\Omega=(a,b)\subseteq\mathbb{R}$ and any $T>0$ there exists $\rho>0$ such that for all 
initial data $u_0,u_1$ satisfying 
$\|u_2\|_{L^2(\Omega)}+\|u_1\|_{H^1(\Omega)}<\rho$,
a solution $u\in C^2(0,T;L^2(\Omega))\cap C^1(0,T;H^1(\Omega))$ to 
\begin{equation}
\begin{split}
& (c^{-2}-2\gamma u)u_{tt}-u_{xx}-\beta u_{txx}=2\gamma (u_t)^2\quad\text{in } (0,T)\times\Omega,\\
& u(t=0)=u_0\,,\ u_t(t=0)=u_1\quad\text{in }\Omega,\\
& u^\beta_n+\sqrt{c^{-2}-2\gamma u}\,u_t
-\frac{\gamma}{2\sqrt{c^{-2}-2\gamma u}}\left(
u_t u -\frac{1}{\sqrt{c^{-2}-2\gamma u}} u^\beta_n u\right)
=0\\
& \hspace*{8cm}\text{at } (0,T)\times\{a,b\}
\end{split}
\label{Westervelt_first_1D_IBVP_intro}
\end{equation}
exists and is unique.
\end{theorem}

\begin{theorem}\label{theo:zero2d}
For $\beta>0$, any smooth and bounded domain $\Omega\subseteq\mathbb{R}^2$ and any $T>0$ there exists $\rho>0$ such that for all 
initial data $u_0,u_1$ satisfying 
$\|u_2\|_{L^2(\Omega)}+\|u_1\|_{H^1(\Omega)}+\|u_0\|_{H^2(\Omega)}<\rho$, 
a solution $u\in C^2(0,T;L^2(\Omega))\cap C^1(0,T;H^1(\Omega))\cap C(0,T;H^2(\Omega))$ to 
\begin{equation}
\begin{split}
& (c^{-2}-2\gamma u)u_{tt}-\Delta u-\beta \Delta u_t=2\gamma (u_t)^2\quad\text{in } (0,T)\times\Omega,\\
& u(t=0)=u_0\,,\ u_t(t=0)=u_1\quad\text{in }\Omega,\\
& u^\beta_n+\sqrt{c^{-2}-2\gamma u}\, u_t =0\quad\text{at } (0,T)\times\partial\Omega
\end{split}
\label{Westervelt_zero_2D_IBVP_intro}
\end{equation}
exists and is unique.
\end{theorem}

\begin{theorem}\label{theo:first2d}
For $\beta>0$, any smooth and bounded domain $\Omega\subseteq\mathbb{R}^2$ and any $T>0$ there exists $\rho>0$ such that for all 
initial data $u_0,u_1$ satisfying 
$\|u_2\|_{L^2(\Omega)}+\|u_1\|_{H^1(\Omega)}+\|u_0\|_{H^2(\Omega)}<\rho$, 
a solution $u\in C^2(0,T;L^2(\Omega))\cap C^1(0,T;H^1(\Omega))\cap C(0,T;H^2(\Omega))$ to 
\begin{equation}
\begin{split}
& (c^{-2}-2\gamma u)u_{tt}-\Delta u-\beta \Delta u_t=2\gamma (u_t)^2\quad\text{in } (0,T)\times\Omega,\\
& u(t=0)=u_0\,,\ u_t(t=0)=u_1\quad\text{in }\Omega,\\
& 
u^\beta_{tn}+\sqrt{c^{-2}-2\gamma u}\frac{u_{tt}+u^\beta_{tt}}{2}
-\frac{1}{2\sqrt{c^{-2}-2\gamma u}} (u^\beta_{\vartheta\vartheta}+\beta u^\beta_{\vartheta\vartheta})\\
&-\frac{\gamma}{2\sqrt{c^{-2}-2\gamma u}}
\left(u_t -\frac{1}{\sqrt{c^{-2}-2\gamma u}} u_n \right)u_{t}
\\
&
-\frac{\gamma}{2(c^{-2}-2\gamma u)^{3/2}}
\left(\frac12 u_t -\frac{1}{\sqrt{c^{-2}-2\gamma u}} u_n \right)\int_0^\cdot (u^\beta_{\vartheta\vartheta}+\beta u^\beta_{\vartheta\vartheta})\, dt
=0\\
&\hspace*{8cm}\text{at } (0,T)\times\partial\Omega
\end{split}
\label{Westervelt_first_2D_IBVP_intro}
\end{equation}
exists and is unique.
\end{theorem}

\medskip

The remainder of this paper is organized as follows. 
In Section \ref{sec:ABCWestervelt} we derive absorbing boundary conditions for the Westervelt equation via (formal) pseudodifferential calculus in one and two space dimensions. 
Section \ref{secwell} is devoted to energy estimates and the proofs of Theorems \ref{theo:zero1d}--\ref{theo:first2d}. 
In Section \ref{sec:numres} we provide numerical results.

\section{Derivation of absorbing boundary conditions for the Westervelt equation in one and two space dimensions}
\label{sec:ABCWestervelt}

In our derivation, without loss of generality we consider the simple domains $\Omega=(-\infty,0]$ in 1-d and 
$\Omega=(-\infty,0]\times\mathbb{R}$ in 2-d, where $x$ plays the role of the outward unit normal and (in 2-d) $y$ is the tangential direction. 
Moreover, we will skip the term $\beta\Delta u_t$ for the reasons outlined in Remark \ref{rem:beta0}.

\subsection{Absorbing boundary conditions in 1-d via linearization and pseudodifferential calculus}
\label{sec:LinWestervelt_1d}
As it was already mentioned, the direct reformulation of~(\ref{Westervelt}) in terms of pseudo-differential operators
is not possible because of the nonlinear term on the right hand side. Therefore we consider some linearization
around a reference solution $u^{(0)}$
\begin{equation}
(c^{-2}-2\gamma u^{(0)}) u_{tt}-\Delta u
= 2\gamma u^{(0)}_t u_t \quad\text{in }(0,T)\times\Omega, 
\label{Westervelt_lin}
\end{equation}
of this equation. After derivation of the ABCs from this inhomogeneous linear wave equation with variable coefficients, we will re-insert $u^{(0)}=u$ to arrive at ABCs for the Westervelt equation.
The reason for using \eqref{Westervelt_lin} (as was also done for the wellposedness proof in \cite{KL09}) and not the standard linearization according to first order Taylor expansion, which would be 
\begin{equation}
c^{-2}u_{tt}-\Delta u
=2\gamma \Bigl(u^{(0)} u_{tt} + u u^{(0)}_{tt} - u^{(0)} u^{(0)}_{tt} + 2 u^{(0)}_t u_t - (u^{(0)}_t)^2\Bigr)\quad\text{in }(0,T)\times\Omega, 
\label{Westervelt_linTaylor}
\end{equation}
is that the offset terms $- 2\gamma u^{(0)} u^{(0)}_{tt} - 2\gamma (u^{(0)}_t)^2=-\gamma(u^{(0)})^2_{tt}$ would lead to problems with the commutativity of the pseudodifferential operators below.

For simplicity of exposition we first of all consider the one-dimensional version of the Westervelt equation~(\ref{Westervelt}) 
\begin{equation}
c^{-2}u_{tt}-u_{xx}
=\gamma(u^2)_{tt}.
\label{Westervelt_1D}
\end{equation}
In 1-d the linearization \eqref{Westervelt_lin} reads as
\begin{equation}
\mathfrak{D}_1u=0,\quad \text{with }\mathfrak{D}_1=\nu^2\partial^2_t-\partial^2_x
- 2\gamma u^{(0)}_t \partial_t,
\label{Westervelt_1D_u1}
\end{equation}
where we  set $\nu^2=\nu^2(u^{(0)})$ with
\begin{equation}\label{nu}
\nu^2(v)=c^{-2}-2\gamma v\,,
\end{equation}
and point out that the analysis of the Westervelt equation is based on estimates that actually make sure positivity of $c^{-2}-2\gamma u$, so that $\nu^2>0$ is a natural assumption.
In order to derive transparent boundary conditions for the linearized Westervelt equation~(\ref{Westervelt_1D_u1}) we make use of the theory of pseudo-differential calculus. 
For the purpose of this formal derivation, $\nu$ is assumed to be a $C^{\infty}$ function both in time and space, as needed for applying pseododofferential calculus. Since we do not prove this smoothess, our derivations are only formal.

The key idea behind the derivation of ABCs is mostly based on the Nirenberg factorization of~(\ref{Westervelt_1D_u1}) written in terms of pseudo-differential operators.
To construct approximate boundary conditions one can factorize the operator $\mathfrak{D}_1$ as
\begin{equation}
\mathfrak{D}_1=-(\partial_x-A)(\partial_x-B)+R,
\label{NFact_Westervelt_1D_u1}
\end{equation}
where $A=A(x,t,D_t)$ and $B=B(x,t,D_t)$ are pseudo-differential operators with symbols $a(x,t,\tau)$ and $b(x,t,\tau)$ from the space
\[S^1=S^1(\mathbb{R}^2)=\left\{f(t,\tau)\in C^{\infty}(\mathbb{R}^2):
\left|\frac{\partial^{\xi}}{\partial t^{\xi}}\frac{\partial^{\sigma}}{\partial\tau^{\sigma}}f(t,\tau)\right|\le C_{\xi,\sigma}(1+|\tau|)^{1-|\sigma|},\ \forall
\xi,\sigma\in\mathbb{N}_0\right\}.\]
The differential operator $D_t$ is defined as $-i\partial_t$ with the imaginary unit $i$, and
$R$ is a smoothing pseudo-differential operator with the Schwartz kernel $k(x,y)\in C^{\infty}$ 
satisfying~\cite{Hormander1985}

\begin{equation*}
(1+|x-y|)^N\left|\frac{\partial^{\xi}}{\partial x^{\xi}}\frac{\partial^{\sigma}}{\partial y^{\sigma}}k(x,y)\right|\le C_{\xi,\sigma,N},\quad 
\forall \xi,\sigma, N\in\mathbb{N}_0.
\end{equation*}

Developing factorization~(\ref{NFact_Westervelt_1D_u1}), we get
\begin{equation}
\mathfrak{D}_1=-\partial^2_x+(A+B)\partial_x+B_x-AB+R.
\label{NirFact_D1_op3}
\end{equation}
At the symbolic level, factorization~(\ref{NirFact_D1_op3}) reduces to
\begin{equation}
\nu^2(i\tau)^2
- 2\gamma u^{(0)}_t (i\tau)
=(a+b)\partial_x+b_x-ab+R
\label{NirFact_D1_sym1}
\end{equation}
with the correspondence $i\tau\leftrightarrow\partial_t$ between the frequency and the (physical) time domains, and 
where by a slight abuse of notation, for a function $f$, we denote the symbol of the zero order differential 
operators $u\mapsto f u$ (multiplication operator) again by $f$. 

Now, we have to define symbols $a$ and $b$ in~(\ref{NirFact_D1_sym1}). For doing so,
it is worth to remark that formally these symbols admit the following asymptotic expansions
\begin{subequations}
\label{AsympExp_ab_1D}
\begin{equation}
a(x,t,\tau)\sim\sum_{j\ge0}a_{1-j}(x,t,\tau),\quad |\tau|\rightarrow\infty
\label{AsympExp_a_1D}
\end{equation}
and
\begin{equation}
b(x,t,\tau)\sim\sum_{j\ge0}b_{1-j}(x,t,\tau),\quad |\tau|\rightarrow\infty \ ,
\label{AsympExp_b_1D}
\end{equation}
\end{subequations}
where $a_{1-j}(x,t,\tau)$ and $b_{1-j}(x,t,\tau)$ are homogeneous of degree $1-j$ in $\tau$. To proceed, one has to
substitute~(\ref{AsympExp_ab_1D}) in~(\ref{NirFact_D1_sym1}) and equate symbols of the same degree of homogeneity on
both sides of equality~(\ref{NirFact_D1_sym1}). However, before this substitution, we recall the reader
the definition of the product of two pseudo-differential operators which are used owing to the term $ab$ in~(\ref{NirFact_D1_sym1}).

In accordance to the theorem on the product of two pseudo-differential operators~\cite{Wong1999},
$A(\mathbf{x},D)\in \Psi^{m_1}$ and $B(\mathbf{x},D)\in \Psi^{m_2}$ with symbols $a(\mathbf{x},\boldsymbol{\xi})\in S^{m_1}$ and $b(\mathbf{x},\boldsymbol{\xi})\in S^{m_2}$ respectively,
a composition operator $C(\mathbf{x},D)=A(\mathbf{x},D)B(\mathbf{x},D)\in\Psi^{m_1+m_2}$ has the asymptotic expansion of its symbol $c(\mathbf{x},\boldsymbol{\xi})\in S^{m_1+m_2}$ given by 
\begin{equation}
c(\mathbf{x},\boldsymbol{\xi})\sim\sum\limits_{|\alpha|\le N}\frac{1}{\alpha!}D^{\alpha}_{\xi}a(\mathbf{x},\boldsymbol{\xi})\partial^{\alpha}_{x}b(\mathbf{x},\boldsymbol{\xi})
\label{PDOproduct1}
\end{equation}
for every nonnegative integer $N$ and with the standard multi-index notation $\alpha=(\alpha_1,\alpha_2,\ldots,\alpha_k)$ and $|\alpha|=\alpha_1+\alpha_2+\ldots+\alpha_k$,
$\mathbf{x}=(x_1,x_2,\ldots,x_k)$, $\boldsymbol{\xi}=(\xi_1,\xi_2,\ldots,\xi_k)$, $D^{\alpha}=D^{\alpha_1}D^{\alpha_2}\ldots D^{\alpha_k}$ and $\partial^{\alpha}=\partial^{\alpha_1}\partial^{\alpha_2}\ldots\partial^{\alpha_k}$.

Thus, the symbol $c:=ab$ of the product of the pseudo-differential operators $A(x,t,D_t)B(x,t,D_t)$, is asymptotic to 
\begin{equation}
c(x,t,\tau)\sim\sum\limits_{k,l,n\ge0}\frac{(-i)^n}{n!}\partial^n_{\tau}a_{1-k}(x,t,\tau)\partial^n_{t}b_{1-l}(x,t,\tau).
\label{PDOproduct2}
\end{equation}

Substitution of~(\ref{AsympExp_ab_1D}) and~(\ref{PDOproduct2}) in~(\ref{NirFact_D1_sym1}) and casting-out $R$ lead to
\begin{equation}
\begin{split}
\nu^2(i\tau)^2
-2\gamma u^{(0)}_t(i\tau)&=\sum\limits_{j\ge0}(a_{1-j}+b_{1-j})\partial_x+\sum\limits_{j\ge0}\partial_xb_{1-j}\\
&-\sum\limits_{j\ge0, \, k+l+n=j}\underbrace{\left(\frac{(-i)^n}{n!}\partial^n_{\tau}a_{1-k}\partial^n_{t}b_{1-l}\right)}_{\mathcal{O}(\tau^{2-j})}
\quad k,l,n\ge0.
\end{split}
\label{NirFact_D1_sym2}
\end{equation}

Evidently, the more coefficients are taken in~(\ref{NirFact_D1_sym2}) the more accurate ABCs are.
However, taking more coefficients also makes the ABCs more complicated and involved to implement since they will contain higher order derivatives.
Therefore, we only show how to find $\{a_j,b_j\}_{j=\{1,0,-1\}}$.
In order to define the first pair of coefficients $a_1$, $b_1$ one has to equate the symbols with the degree of homogeneity $\mathcal{O}(\tau^2)$. 
This gives the following system of equations
\begin{equation}
 \left\{
  \begin{aligned}
   a_1+b_1&=0\\
   \nu^2(i\tau)^2&=-a_1b_1.
  \end{aligned}
  \right.
\label{a1b1_eq} 
\end{equation}
The solutions to~(\ref{a1b1_eq}) are given by
\begin{equation}
b^{(1,2)}_1=-a^{(1,2)}_1=\pm \nu(i\tau).
\label{a1b1_1} 
\end{equation}
We take
\begin{equation}
b_1=-a_1=\nu(i\tau)
\label{a1b1_1D_1} 
\end{equation}
to make the terms of order $\mathcal{O}(\tau^2)$ vanish.
\begin{remark}
The choice of the sign in front of $\nu(i\tau)$ is not arbitrary.
This sign defines the propagation direction of the wave. 
\end{remark}

In order to find the next pair of coefficients $a_0$, $b_0$ we equate symbols with degree of homogeneity $\mathcal{O}(\tau^1)$. 
In other words, we have to solve the system
\begin{equation}
 \left\{
  \begin{aligned}
   a_0+b_0&=0,\\
	2\gamma u^{(0)}_t(i\tau)&=a_1b_0+a_0b_1-i{a_1}_{\tau}{b_1}_t-{b_1}_x,
  \end{aligned}
  \right.
\label{a0b0_1D_tau1}
\end{equation}
in terms of unknown $a_0$, $b_0$. 
Substitution of $b_1=-a_1$ in~(\ref{a0b0_1D_tau1}) gives
\begin{equation}
b_0=-a_0=-\frac{1}{2a_1}\left(i{a_1}_\tau {a_1}_t+{a_1}_x
-2\gamma u^{(0)}_t(i\tau)\right)
\label{a0b0_1D_eq2_3}
\end{equation}
or in terms of $a_1=-\nu(i\tau)$ we have
\begin{equation}
b_0=-a_0=-\frac{1}{2\nu}\left(\mathcal{A}_0[\nu]
+2\gamma u^{(0)}_t\right)
\label{a0b0_1D_eq2_5}
\end{equation}
with the operator $\mathcal{A}_0:=\partial_x+\nu\partial_t$. 

In order to obtain more accurate boundary conditions one has to equate the symbols with degree of homogeneity $\mathcal{O}(\tau^0)$
which leads to the following system
\begin{equation}
 \left\{
  \begin{aligned}
   a_{-1}+b_{-1}&=0,\\
   -a_1b_{-1}-a_0b_0-a_{-1}b_1+i({a_1}_{\tau}{b_0}_t+{a_0}_{\tau}{b_1}_t)-\frac{i^2}{2}{a_1}_{\tau\tau}{b_1}_{tt}+{b_0}_x&=0.
  \end{aligned}
  \right.
\label{am1bm1_1D_tau0_1}
\end{equation}
The solution of \eqref{am1bm1_1D_tau0_1} is given by
\begin{equation}
b_{-1}=-a_{-1}=-\frac{1}{2a_1}\left(-a_0^2+i({a_1}_{\tau}{a_0}_t+{a_0}_{\tau}{a_1}_t)+\frac{1}{2}{a_1}_{\tau\tau}{a_1}_{tt}+{a_0}_x
\right).
\label{am1bm1_1D_tau0_2}
\end{equation}
Taking into account~(\ref{a1b1_1D_1}) and~(\ref{a0b0_1D_eq2_5}) we deduce that
\begin{equation}
\begin{split}
&b_{-1}=-a_{-1}\\
&=\frac{1}{2\nu(i\tau)}
\left(\mathcal{A}_0\left[\frac{1}{2\nu}\left(\mathcal{A}_0[\nu]
+2\gamma u^{(0)}_t\right)\right]\right.
\left. -\left(\frac{1}{2\nu}\left(\mathcal{A}_0[\nu]
+2\gamma u^{(0)}_t\right)\right)^2
\right)\\
&=:\frac{\gamma\mu}{2\nu(i\tau)}\ .
\end{split}
\label{am1bm1_1D_tau0_3}
\end{equation}
Note that with the Taylor linearization \eqref{Westervelt_linTaylor} an offset term $\gamma\underline{({u^{(0)}}^2)_{tt}}$ would have appeared here which would have prevented the equality $a_{-1}a_1=a_1a_{-1}$.
(Here, we write $\underline{f}$ for the symbol of the zero order differential operator $u\mapsto f $ (constant mapping), which has to be strictly distinguished from the multiplication operator $u\mapsto f u$.)
This problem is avoided by using the fixed point type linearization \eqref{Westervelt_lin}.

In accordance to~\cite{MajdaOsher1975},
the operator
\begin{equation}
\partial_x-a(x,t,D_t)=0
\label{EM_ABC_1D_op}
\end{equation}
annihilates outgoing waves at $\{x=0\}\times(0,T)$. 
Substitution of the asymptotic expansion~(\ref{AsympExp_a_1D}) with the first $k$ leading terms results in the following boundary condition
\begin{equation}
\left.\left(\partial_x-\sum\limits^k_{j=0}a_{1-j}(x,t,D_t)\right)u\right|_{x=0}=0,
\label{ABC_1D_sym1}
\end{equation}
i.e., an ABC of order $k$ is obtained by keeping the first $k$ terms in the asymptotic expansions~(\ref{AsympExp_ab_1D}).
 
Thus in order to construct a zero order ABC we set $k=0$ and substitute the coefficient $a_1$ in~(\ref{ABC_1D_sym1}) which gives
%
\label{ABC0_1D_x0xa}
\begin{equation}
\mathcal{A}_0[u]\Big|_{x=0}=\left(u_x+\nu u_t\right)\Big|_{x=0} =0.
\label{ABC0_1D_xa}
\end{equation}

Parallel to the construction of the zero order ABCs~(\ref{ABC0_1D_x0xa}), we set $k=1$ and substitute $a_1$, $a_0$
in~(\ref{ABC_1D_sym1}) to obtain the first order boundary conditions:
\label{ABC1_1D_x0xa}
\begin{equation}
\mathcal{A}_1u\Big|_{x=0}=(\mathcal{A}_0-\mathcal{B}_1)u\Big|_{x=0}=
\left.\left(u_x+\nu u_t-\frac{1}{2\nu}\left((\nu_x+\nu \nu_t)u+2\gamma u^{(0)}_t u\right)\right)
\right|_{x=0}=0
\label{ABC1_1D_xa}
\end{equation}
with $\mathcal{B}_1:=\frac{1}{2\nu}\left(\mathcal{A}_0[\nu]
+2\gamma u^{(0)}_t\right)$.

For $k=2$ we obtain the second order ABCs
\begin{equation}
\begin{aligned}
&\mathcal{A}_2u\Big|_{x=0}=(\mathcal{A}_1u_t-\mathcal{B}_2 u)\Big|_{x=0}\\
&=
\left.\left(u_{xt}+\nu u_{tt}-\frac{1}{2\nu}\left((\nu_x+\nu \nu_t)u_t+2\gamma u^{(0)}_t u_t-\mu u\right)\right)
\right|_{x=0}=0,
\label{ABC2_1D_xa}
\end{aligned}
\end{equation}
where we have multiplied with $(i\tau)$ before converting from symbols to operators, and where
$\mathcal{B}_2:=\frac{\gamma\mu(u^{(0)})}{2\nu(u^{(0)})}$ with
\begin{eqnarray}\label{mu}
\mu(v)&=&\frac{1}{\gamma}
\mathcal{A}_0\left[\frac{1}{2\nu(v)}\left(\mathcal{A}_0[\nu(v)]
+2\gamma v_t\right)\right]
 -\left(\frac{1}{2\nu(v)}\left(\mathcal{A}_0[\nu(v)]
+2\gamma v_t\right)\right)^2\nonumber\\
&=&
\mathcal{A}_0\left[\frac{1}{2\sqrt{c^{-2}-2\gamma v}}
\left(-\frac{v_x}{\sqrt{c^{-2}-2\gamma v}}+ v_t\right)\right]
\nonumber\\
&& -\gamma\left(\frac{1}{2\sqrt{c^{-2}-2\gamma v}}
\left(-\frac{ v_x}{\sqrt{c^{-2}-2\gamma v}}+ v_t\right)\right)^2\,.
\end{eqnarray}

Inserting $u$ itself for the a priori solution $u^{(0)}$, we arrive at zero
\begin{equation}
\left.\left(u_x+\sqrt{c^{-2}-2\gamma u}\ u_t\right)\right|_{x=0} =0,
\label{ABC0_1D}
\end{equation}
first
\begin{equation}
\left.\left(u_x+\sqrt{c^{-2}-2\gamma u}\ u_t-\frac{\gamma}{2\sqrt{c^{-2}-2\gamma u}}\left(
u_t u -\frac{1}{\sqrt{c^{-2}-2\gamma u}} u_x u\right)\right)
\right|_{x=0}=0 
\label{ABC1_1D}
\end{equation}
and second order
\begin{equation}
\begin{aligned}
&\Bigl(u_{xt}+\sqrt{c^{-2}-2\gamma u}\, u_{tt}\\
&\left.-\frac{\gamma}{2\sqrt{c^{-2}-2\gamma u}}\Bigl(
(u_t)^2 -\frac{1}{\sqrt{c^{-2}-2\gamma u}}\, u_x u_t-\mu(u) u\Bigr)\Bigr)
\right|_{x=0}=0 
\label{ABC2_1D}
\end{aligned}
\end{equation}
nonlinear ABCs.
We will see in Section \ref{sec_paradiff} that slightly different conditions result from derivation via a paradifferential approach.

\subsection{Absorbing boundary conditions in 2-d via linearization and pseudodifferential calculus}
\label{sec:LinWestervelt_2d}
In the spatially two dimensional situation 
\begin{equation}
\mathfrak{D}_1u=0,\quad \text{with }\mathfrak{D}_1
=\nu^2\partial^2_t-\partial^2_x-\partial^2_y
- 2\gamma u^{(0)}_t \partial_t,
\label{Westervelt_2D_u1}
\end{equation}
on the domain $(-\infty,0)\times\mathbb{R}$, where $\nu$ is defined by \eqref{nu}, we proceed very similarly to the 1-d case:
We consider pseudo-differential operators $A=A(x,y,t,D_y,D_t)$ and $B=B(x,y,t,D_y,D_t)$
with respect to time and tangential (i.e., $y$) direction, but the expansion is still with respect to powers of $\tau$, so equations 
\eqref{NFact_Westervelt_1D_u1}, \eqref{NirFact_D1_op3}
(with $A=A(x,y,t,D_y,D_t)$ and $B=B(x,y,t,D_y,D_t)$) remain the same whereas 
\eqref{NirFact_D1_sym1}, \eqref{AsympExp_ab_1D}, \eqref{NirFact_D1_sym2} change to
\begin{equation}
\nu^2(i\tau)^2-(i\eta)^2
- 2\gamma u^{(0)}_t (i\tau)
=(a+b)\partial_x+b_x-ab+R
\label{NirFact_D1_sym1_2d}
\end{equation}
with the correspondence $i\eta\leftrightarrow\partial_y$ and
\begin{subequations}
\label{AsympExp_ab_2D}
\begin{equation}
a(x,y,t,\eta,\tau)\sim\sum_{j\ge0}a_{1-j}(x,y,t,\eta,\tau),\quad |\tau|\rightarrow\infty
\label{AsympExp_a_2D}
\end{equation}
\begin{equation}
b(x,y,t,\eta,\tau)\sim\sum_{j\ge0}b_{1-j}(x,y,t,\eta,\tau),\quad |\tau|\rightarrow\infty
\label{AsympExp_b_2D}
\end{equation}
\end{subequations}
and
\begin{equation}
\begin{split}
\nu^2(i\tau)^2-(i\eta)^2
-2\gamma u^{(0)}_t(i\tau)
&=\sum\limits_{j\ge0}(a_{1-j}+b_{1-j})\partial_x+\sum\limits_{j\ge0}\partial_xb_{1-j}\\
&-\sum\limits_{j\ge0, \, k+l+n=j}\underbrace{\left(\frac{(-i)^n}{n!}\partial^n_{\tau}a_{1-k}\partial^n_{t}b_{1-l}\right)}_{\mathcal{O}(\tau^{2-j})}
\quad k,l,n\ge0,
\end{split}
\label{NirFact_D1_sym2_2d}
\end{equation}
respectively, where $a_{1-j}$ and $b_{1-j}$ are homogeneous of degree $1-j$ in $\tau$ (and are additionally functions of $x,y,t$, and $\eta$).
As in \cite{EngquistMajda1977}, in our derivations we will rely on an assumption of the type $\eta\sim\tau$ or even $\frac{\eta}{\tau}$ small. 
Considering the $\mathcal{O}(\tau^2)$ terms in \eqref{NirFact_D1_sym2_2d} leads us to   
\begin{equation}
 \left\{
  \begin{aligned}
   \nu^2(i\tau)^2-(i\eta)^2&=-a_1b_1,\\
   a_1+b_1&=0.
  \end{aligned}
  \right.
\label{a1b1_eq_2d} 
\end{equation}
in place of \eqref{a1b1_eq}, which leads to 
\begin{equation}
b_1=-a_1=\sqrt{\nu^2(i\tau)^2-(i\eta)^2}.
\label{a1b1_2D_1} 
\end{equation}
At this point, a fundamental difference to the 1-d case arises, since we will have to approximate the square root 
\[\sqrt{\nu^2(i\tau)^2-(i\eta)^2}=\nu(i\tau)\sqrt{1-\frac{\eta^2}{\nu^2\tau^2}}\]
in order to derive practically applicable boundary conditions. We will do so by a Taylor expansion whose order is adapted to the order of the ABCs.

The computations for $a_0$, $b_0$ look exactly the same as in the 1-d case and yield
\begin{equation}
b_0=-a_0=-\frac{1}{2a_1}\left(i{a_1}_\tau {a_1}_t+{a_1}_x
-2\gamma u^{(0)}_t(i\tau)\right)
\label{a0b0_2D_eq2_3}
\end{equation}
i.e., 
\begin{equation}
b_0=-a_0=-\frac{\nu_t}{2}\left(1-\frac{\eta^2}{\nu^2\tau^2}\right)^{-3/2}
- \frac{\nu_x}{2\nu}\left(1-\frac{\eta^2}{\nu^2\tau^2}\right)^{-1}
-\frac{
2\gamma u^{(0)}_t}{2\nu}\left(1-\frac{\eta^2}{\nu^2\tau^2}\right)^{-1/2}.
\label{a0b0_2D_eq2_4}
\end{equation}

To obtain zero order boundary conditions we use the zero order Taylor expansion 
\[(1-x)^{1/2} \approx 1,\quad x:=\frac{\eta^2}{\nu^2\tau^2} \]
in \eqref{a1b1_2D_1} to end up with 
\begin{equation}
\tilde{b}_1^0=-\tilde{a}_1^0=\nu(i\tau).
\label{tilde_a0b0_2D}
\end{equation}
For our first order boundary conditions we use the first order Taylor approximations 
\[(1-x)^{1/2} \approx 1- \frac12 x \,, \quad 
(1-x)^{-3/2} \approx 1+\frac32 x \,, \quad 
(1-x)^{-1} \approx 1+ x \,, \quad 
(1-x)^{-1/2} \approx 1+ \frac12 x 
\]
for the terms that are nonlinear with respect to $\tau$, $\eta$ in \eqref{a0b0_2D_eq2_3}, \eqref{a0b0_2D_eq2_4}. This yields the symbols
\begin{eqnarray*}
\tilde{b}_1^1=-\tilde{a}_1^1&=&\nu(i\tau)\left(1-\frac{\eta^2}{2\nu^2\tau^2}\right),\\
\tilde{b}_0^1=-\tilde{a}_0^1&=&
-\frac{\nu_t}{2}\left(1+\frac{3\eta^2}{2\nu^2\tau^2}\right)
- \frac{\nu_x}{2\nu}\left(1+\frac{\eta^2}{\nu^2\tau^2}\right)
-\frac{2\gamma u^{(0)}_t}{2\nu}\left(1+\frac{\eta^2}{2\nu^2\tau^2}\right).
\end{eqnarray*}

Again we insert $u$ itself for the a priori solution $u^{(0)}$ to arrive at zero order ABCs
\begin{equation}
\left.\left(u_x+\sqrt{c^{-2}-2\gamma u}\, u_t\right)\right|_{x=0} =0
\label{ABC0_2D}
\end{equation}
and at first order ABCs
\begin{equation}
\begin{aligned}  
&\left((u_{xt}+\sqrt{c^{-2}-2\gamma u}\, u_{tt}-\frac{1}{2\sqrt{c^{-2}-2\gamma u}} u_{yy}
\right.\\
&-\frac{\gamma}{2\sqrt{c^{-2}-2\gamma u}}
\left(u_t -\frac{1}{\sqrt{c^{-2}-2\gamma u}} u_x \right)u_{t}
\\
&\left.\left.+\frac{\gamma}{2(c^{-2}-2\gamma u)^{3/2}}
\left(\frac12 u_t +\frac{1}{\sqrt{c^{-2}-2\gamma u}} u_x \right) \int_0^\cdot u_{yy} \, dt\right)
\right|_{x=0}=0,
\end{aligned} 
\label{ABC1_2D}
\end{equation}
where we have multiplied the symbols with $(i\tau)$ to obtain \eqref{ABC1_2D}. 

\section{Well-posedness}\label{secwell}
In this section we will show well-posedness of the Westervelt equation with zero or first order ABC derived above in one or two space dimensions.
Note that zero order ABC have already been considered in \cite{Clason2009}. However, the conditions there do not take into account the nonlinearity in 
the highest order time derivative. Moreover, the proof in \cite{Clason2009} is carried out in higher space dimensions, which necessitates the use of higher 
order energies. In 1-d this is not required (simply due to the fact that in 1-d already $H^1$ embeds into $L^\infty$) and the proof is on one hand much simpler, on the other hand it enables existence also of spatially less smooth solutions and well-posedness in the absence of interior damping (i.e., with $\beta=0$). 
For these reason we will also provide the well-posedness proof for the 1-d Westervelt equation with zero order ABC \eqref{ABC0_1D} here.

Since we derive energy estimates by only multiplying with $u_t$ for zero order ABC in 1-d, the strong damping term and the terms resulting from its integration by parts at the boundary will be easily tractable in that case.
However, for the first order ABCs, carrying out energy estimates following the idea in \cite{HJ994}, the $\beta$ term yields derivatives of $u$ on the boundary that are too high to be controllable by other boundary (or, via trace theorems, interior) terms. Therefore we will modify the first order ABC accordingly to account for the strong damping and arrive at decaying energies. Note that in the derivations of section \ref{sec:ABCWestervelt} we had omitted the $\beta$ terms since they would have destroyed commutativity. The terms that we insert now again in favor of energy decay are different from those omitted in section \ref{sec:ABCWestervelt}, though. I.e., the (formal) Nirenberg factorization from there would not have helped in obtaining energy dissipation. In fact it turns out that the ABCs derived in section \ref{sec:ABCWestervelt} (plus the $\beta$-modifications made here) only allow us to show local in time well-posedness. 
As to be expected, the resulting ABCs coincide with the classical Engquist-Majda ones in case of constant coefficients and vanishing damping.

\subsection{Energy identities for the strongly damped inhomogeneous wave equation with variable coefficients} \label{subsec:beta}
Before proceeding to well-posedness of the nonlinear Westervelt equation with zero and first order ABCs, we will derive some energy identities (especially we will carry over the energy identities used in the well-posedness proof for first order ABCs in \cite{HJ994}) for inhomogeneous wave equations with variable coefficients and strong damping of the form
\begin{equation} \label{PDEalpha}
\alpha u_{tt} -\Delta u -\beta\Delta u_t =f u_t +g\quad \text{ in }(0,T)\times \Omega
\end{equation}
with 
$\Omega\subseteq\mathbb{R}^{d}$, 
$d\in\{1,2\}$ 
$\alpha=\alpha(t,x(,y))>0$, $f=f(t,x(,y))$, 
$\beta\geq0$ and initial conditions $u(t=0)=u_0$, $u_t(t=0)=u_1$.
This will provide us with crucial information on how to incorporate the strong interior damping term into the ABCs and help 
to prove well-posedness of the Westervelt equation with ABCs in the next subsections.
For simplicity of exposition we will here restrict ourselves to a geometry with $x$ being the boundary normal and $y$ the boundary tangential direction, respectively. The general case can be covered by applying smooth local boundary transformations.

Multiplying the PDE with $u_t$ we obtain
\begin{eqnarray}\label{energyid0}
\lefteqn{
\frac12 \Bigl(\|\sqrt{\alpha(t)} u_t(t)\|_{L^2(\Omega)}^2 + \|\nabla u(t)\|_{L^2(\Omega)}^2 \Bigr)
+\beta\int_0^t \|\nabla u_t\|_{L^2(\Omega)}^2\, ds} \nonumber\\
&=& 
\frac12 \Bigl(\|\sqrt{\alpha(0)} u_t(0)\|_{L^2(\Omega)}^2 + \|\nabla u(0)\|_{L^2(\Omega)}^2\Bigr)\nonumber\\
&&+\int_0^t\int_\Omega \Bigl( \alpha u_{tt} u_t + \frac12\alpha_t (u_t)^2
+ \nabla u_t \nabla u + \beta |\nabla u_t|^2 \nonumber\\
&=& 
\frac12 \Bigl(\|\sqrt{\alpha(0)} u_t(0)\|_{L^2(\Omega)}^2 + \|\nabla u(0)\|_{L^2(\Omega)}^2\Bigr)\nonumber\\
&&+\int_0^t\int_\Omega \Bigl(\Bigl( \alpha u_{tt}  -\Delta u -\beta\Delta u_t\Bigr) u_t 
+ \frac12\alpha_t (u_t)^2 \Bigr)ds
+\int_0^t\int_{\partial\Omega} \Bigl(u+\beta u_t\Bigr)_n u_t \, d\Gamma\, dt \nonumber\\
&=& 
\frac12 \Bigl(\|\sqrt{\alpha(0)} u_t(0)\|_{L^2(\Omega)}^2 + \|\nabla u(0)\|_{L^2(\Omega)}^2\Bigr)\nonumber\\
&&+\int_0^t\int_\Omega \Bigl\{\Bigl(f + \frac12\alpha_t\Bigr) (u_t)^2 +g u_t\Bigr\} ds
+\int_0^t\int_{\partial\Omega} u^\beta_n u_t \, d\Gamma\, dt \,,
\end{eqnarray}
where 
$$u^\beta=u+\beta u_t\,.$$

This suggests to use as zero order absorbing boundary conditions
\begin{equation}\label{0orderABC_energy}
\sqrt{\alpha} u_t+u^\beta_n = 0 \mbox{ (or }\sqrt{\alpha} u_t+u^\beta_n =\mbox{ lower order terms),} 
\end{equation}
where  ``lower order terms'' are expressions whose $L^2(0,T;L^2(\partial\Omega))$ inner product with $u_t$ can be dominated by the energy 
$$
E_0[u](t)=
\frac12 \Bigl(\|\sqrt{\alpha(t)} u_t(t)\|_{L^2(\Omega)}^2 + \|\nabla u(t)\|_{L^2(\Omega)}^2 \Bigr)
$$
and/or the interior dissipation
$$
\beta\int_0^t \|\nabla u_t\|_{L^2(\Omega)}^2\, ds\,,
$$
and/or the boundary dissipation
$$
\int_0^t \|\alpha^{-1/4}u^\beta_n\|_{L^2(\partial\Omega)}^2\, ds
=\int_0^t \|\alpha^{1/4}u_t\|_{L^2(\partial\Omega)}^2\, ds\,.
$$

Similarly, if we differentiate the PDE wrt $t$ and multiply with $u_{tt}$, we arrive (after space and time integration) at the energy identity
\begin{eqnarray}\label{energyid1}
\lefteqn{
E_1[u](t)
+\beta\int_0^t \|\nabla u_{tt}\|_{L^2(\Omega)}^2\, ds} \nonumber\\
&=& 
\frac12 \Bigl(\|\sqrt{\alpha(0)} u_{tt}(0)\|_{L^2(\Omega)}^2 + \|\nabla u_t(0)\|_{L^2(\Omega)}^2\Bigr)\nonumber\\
&&+\int_0^t\int_\Omega \Bigl( \alpha u_{ttt} u_{tt} + \frac12\alpha_t (u_{tt})^2
+ \nabla u_{tt} \nabla u_t + \beta |\nabla u_{tt}|^2 \nonumber\\
&=& 
E_1[u](0)+\int_0^t\int_\Omega \Bigl(\Bigl( (\alpha u_{tt})_t  -\Delta u_t -\beta\Delta u_{tt}\Bigr) u_{tt} 
- \frac12\alpha_t (u_{tt})^2 \Bigr)ds
\nonumber\\
&&
+\int_0^t\int_{\partial\Omega} \Bigl(u_t+\beta u_{tt}\Bigr)_n u_{tt} \, d\Gamma\, dt \nonumber\\
&=& 
E_1[u](0)+\int_0^t\int_\Omega \Bigl\{\Bigl(f - \frac12\alpha_t\Bigr) (u_{tt})^2 + (f_t u_t + g_t) u_{tt}\Bigr\} ds
\nonumber\\
&&
+\int_0^t\int_{\partial\Omega} u^\beta_{tn} u_{tt} \, d\Gamma\, dt \,,
\end{eqnarray}
where 
\begin{equation}\label{E1}
E_1[u](t)=E_0[u_t](t)=\frac12 \Bigl(\|\sqrt{\alpha(t)} u_{tt}(t)\|_{L^2(\Omega)}^2 + \|\nabla u_t(t)\|_{L^2(\Omega)}^2 \Bigr)\,.
\end{equation}

Multiplication of the time differentiated PDE with $\alpha u_{tt}$ (instead of $u_{tt}$) yields
\begin{eqnarray}\label{energyid2}
\lefteqn{
\frac12 \Bigl(\|\alpha(t) u_{tt}(t)\|_{L^2(\Omega)}^2 + \|\sqrt{\alpha(t)} \nabla u_t(t)\|_{L^2(\Omega)}^2 \Bigr)
+\beta\int_0^t \|\sqrt{\alpha} \nabla u_{tt}\|_{L^2(\Omega)}^2\, ds}\nonumber\\
&=& 
\frac12 \Bigl(\|\alpha(0) u_{tt}(0)\|_{L^2(\Omega)}^2 + \|\sqrt{\alpha(0)} \nabla u_t(0)\|_{L^2(\Omega)}^2\Bigr)\nonumber\\
&&+\int_0^t\int_\Omega \Bigl( (\alpha u_{tt})_t \alpha u_{tt} + \frac12\alpha_t |\nabla u_t|^2
+ \underbrace{\alpha \nabla u_{tt} \nabla u_t + u_{tt} \nabla \alpha \nabla u_t}_{=\nabla(\alpha u_{tt})\nabla u_t}
 - u_{tt} \nabla \alpha \nabla u_t
\nonumber\\
&&\qquad
+ \beta
\underbrace{(\alpha |\nabla u_{tt}|^2 + u_{tt} \nabla \alpha \nabla u_{tt})}_{\nabla(\alpha u_{tt})\nabla u_{tt}}
-\beta u_{tt} \nabla \alpha \nabla u_{tt}\Bigr)\, d\Omega\, ds 
\nonumber\\
&=& 
\frac12 \Bigl(\|\alpha(0) u_{tt}(0)\|_{L^2(\Omega)}^2 + \|\sqrt{\alpha(0)} \nabla u_t(0)\|_{L^2(\Omega)}^2\Bigr)
\nonumber\\
&&+\int_0^t\int_\Omega \Bigl( (\alpha u_{tt})_t - \Delta u_t - \beta\Delta u_{tt}\Bigr)\alpha u_{tt} \, d\Omega\, ds
\nonumber\\
&&+\int_0^t\int_\Omega \Bigl( 
\frac12\alpha_t |\nabla u_t|^2 - u_{tt} \nabla \alpha \nabla u_t - \beta u_{tt} \nabla \alpha \nabla u_{tt}
\Bigr) \, d\Omega\, ds
\nonumber\\
&&+\int_0^t\int_{\partial\Omega} \Bigl(u_t+\beta u_{tt}\Bigr)_n \alpha u_{tt} \, d\Gamma\, dt 
\nonumber\\
&=& 
\frac12 \Bigl(\|\alpha(0) u_{tt}(0)\|_{L^2(\Omega)}^2 + \|\sqrt{\alpha(0)} \nabla u_t(0)\|_{L^2(\Omega)}^2\Bigr)
\nonumber\\
&&+\int_0^t\int_\Omega \Bigl( \alpha f_t u_t u_{tt} +\alpha f (u_{tt})^2 + \alpha g_t u_{tt}
+\frac12\alpha_t |\nabla u_t|^2 - u_{tt} \nabla \alpha \nabla u^\beta_t 
\Bigr) \, d\Omega\, ds
\nonumber\\
&&+\int_0^t\int_{\partial\Omega} u^\beta_{tn} \alpha u_{tt} \, d\Gamma\, dt 
\end{eqnarray}

Considering the PDE that results from \eqref{PDEalpha} for $u^\beta$, 
\begin{equation}\label{PDEubeta}
\alpha u^\beta_{tt} -\Delta u^\beta -\beta\Delta u^\beta_t 
=f u^\beta_t +g+\beta g_t+\beta f_t u_t-\beta\alpha_t u_{tt}
\end{equation}
differentiating  wrt $x$
\begin{equation}\label{PDEubeta_x}
\begin{aligned}
&\alpha u^\beta_{ttx} -\Delta u^\beta_x -\beta\Delta u^\beta_{tx} \\
&=f u^\beta_{tx}+f_x u^\beta_t +g_x+\beta g_{tx}+\beta f_t u_{tx}+\beta f_{tx} u_t\\
&\quad-\alpha_t (u^\beta_{tx}-u_{tx})-\beta\alpha_{tx} u_{tt}-\alpha_x u^\beta_{tt}
\quad \text{ in }(0,T)\times \Omega\,,
\end{aligned}
\end{equation}
and multiplying with $u^\beta_{tx}$, we get the energy identity
\begin{eqnarray}\label{energyid3}
\lefteqn{
\frac12 \Bigl(\|\sqrt{\alpha(t)} u^\beta_{tx}(t)\|_{L^2(\Omega)}^2 + \|\nabla u^\beta_x(t)\|_{L^2(\Omega)}^2 \Bigr)
+\beta\int_0^t \|\nabla u^\beta_{tx}\|_{L^2(\Omega)}^2\, ds}
\nonumber\\
&=& 
\frac12 \Bigl(\|\sqrt{\alpha(0)} u^\beta_{tx}(0)\|_{L^2(\Omega)}^2 + \|\nabla u^\beta_x(0)\|_{L^2(\Omega)}^2\Bigr)
\nonumber\\
&&+\int_0^t\int_\Omega \Bigl( \alpha u^\beta_{ttx} u^\beta_{tx} +\frac12 \alpha_t (u^\beta_{tx})^2 
+ \nabla u^\beta_{tx} \nabla u^\beta_x 
+ \beta |\nabla u^\beta_{tx}|^2 \Bigr)\, d\Omega\, ds 
\nonumber\\
&=& 
\frac12 \Bigl(\|\sqrt{\alpha(0)} u^\beta_{tx}(0)\|_{L^2(\Omega)}^2 + \|\nabla u^\beta_x(0)\|_{L^2(\Omega)}^2\Bigr)
\nonumber\\
&=& 
+\int_0^t\int_\Omega \Bigl( \alpha u^\beta_{ttx} - \Delta u^\beta_x - \beta\Delta u^\beta_{tx}
+\frac12 \alpha_t u^\beta_{tx}\Bigr) u^\beta_{tx} \, d\Omega\, ds
\nonumber\\
&&+\int_0^t\int_{\partial\Omega} \Bigl(u^\beta_x+\beta u^\beta_{tx}\Bigr)_n u^\beta_{tx} \, d\Gamma\, dt 
\nonumber\\
&=& 
\frac12 \Bigl(\|\alpha(0) u^\beta_{tt}(0)\|_{L^2(\Omega)}^2 + \|\sqrt{\alpha(0)} \nabla u^\beta_t(0)\|_{L^2(\Omega)}^2\Bigr)
\nonumber\\
&&+\int_0^t\int_\Omega \Bigl( 
(f-\tfrac{\alpha_t}{2}) (u^\beta_{tx})^2+(\beta f_t +\alpha_t)u_{tx}u^\beta_{tx}
-\beta\alpha_{tx} u_{tt}u^\beta_{tx}-\alpha_x u^\beta_{tt}u^\beta_{tx}
\nonumber\\
&& \qquad\qquad+f_x u^\beta_t u^\beta_{tx} 
+\beta f_{tx} u_t u^\beta_{tx}
+(g_x+\beta g_{tx})u^\beta_{tx}\Bigr) \, d\Omega\, ds
\nonumber\\
&&+\int_0^t\int_{\partial\Omega} \Bigl(u^\beta_x+\beta u^\beta_{tx}\Bigr)_n u^\beta_{tx} \, d\Gamma\, ds 
\end{eqnarray}
For the combined higher order energy functional
\[
\begin{aligned}
E_2[u](t)= 
\frac12 \Bigl(&\|\alpha(t) u_{tt}(t)\|_{L^2(\Omega)}^2 + \|\sqrt{\alpha(t)} \nabla u_t(t)\|_{L^2(\Omega)}^2\\
&+ \|\sqrt{\alpha(t)} u^\beta_{tx}(t)\|_{L^2(\Omega)}^2 + \|\nabla u^\beta_x(t)\|_{L^2(\Omega)}^2 \Bigr)
\end{aligned}
\]
the identities \eqref{energyid2}, \eqref{energyid3} yield
\begin{eqnarray}\label{energyid23}
\lefteqn{E_2[u](t)
+\beta\int_0^t \Bigl(\|\sqrt{\alpha} \nabla u_{tt}\|_{L^2(\Omega)}^2+\|\nabla u^\beta_{tx}\|_{L^2(\Omega)}^2\Bigr)\, ds 
\ = \ E_2[u](0)
}\nonumber\\
&&+\int_0^t \int_\Omega \Bigl( 
\alpha f (u_{tt})^2 
+\frac12\alpha_t |\nabla u_t|^2 
+(f-\tfrac{\alpha_t}{2}) (u^\beta_{tx})^2
+(\beta f_t +\alpha_t)u_{tx}u^\beta_{tx}
\nonumber\\
&& \qquad\qquad
-\beta\alpha_{tx} u_{tt}u^\beta_{tx}
-\alpha_x u^\beta_{tt}u^\beta_{tx}
+\alpha f_t u_t u_{tt} 
- u_{tt} \nabla \alpha \nabla u^\beta_t 
\nonumber\\
&& \qquad\qquad
+f_x u^\beta_t u^\beta_{tx} 
+\beta f_{tx} u_t u^\beta_{tx}
+ \alpha g_t u_{tt}+(g_x+\beta g_{tx})u^\beta_{tx}
\Bigr) \, d\Omega\, ds
\nonumber\\
&&+\int_0^t\int_{\partial\Omega} \Bigl(\alpha u_{tt}+u^\beta_{xx}+\beta u^\beta_{txx}\Bigr) u^\beta_{tn} \, d\Gamma\, ds \,,
\end{eqnarray}
where we have used the fact that $x$ is the outward normal direction in our setting.
This suggests to use first order boundary conditions leading to the identity 
\begin{equation}\label{bndyid}
\alpha u_{tt}+u^\beta_{xx}+\beta u^\beta_{txx} + 2\sqrt{\alpha}u^\beta_{tn} =0
\mbox{ (or lower order terms)} 
\end{equation}
where this time ``lower order terms'' are expressions whose $L^2(0,T;L^2(\partial\Omega))$ inner product with $u^\beta_{tn}$ can be dominated by the higher order energy $E_2[u](t)$ 
and/or the interior dissipation
$$
\beta\int_0^t \Bigl(\|\sqrt{\alpha} \nabla u_{tt}\|_{L^2(\Omega)}^2+\|\nabla u^\beta_{tx}\|_{L^2(\Omega)}^2\Bigr)\, ds\,,
$$
and/or the boundary dissipation
$$
\int_0^t \|\alpha^{1/4} u^\beta_{tn}\|_{L^2(\partial\Omega)}^2\, ds\,.
$$
Using the PDE \eqref{PDEubeta} for transforming second order normal (i.e., $x$) derivatives to tangential (i.e., $y$) derivatives, we can achieve the boundary identity \eqref{bndyid} e.g. by the first order ABCs
\begin{equation}\label{1orderABC_energy}
\begin{aligned}
&\alpha (u+u^\beta)_{tt}  - u^\beta_{yy}-\beta u^\beta_{tyy}  
-fu^\beta_t-\beta f_tu_t
-g-\beta g_t
+\beta\alpha_t u_{tt} +2\sqrt{\alpha}u^\beta_{tn}=0 \\
&\mbox{ (or lower order terms),}
\end{aligned}\end{equation}
where in one space dimension the $yy$ derivative terms are just skipped.

In the nonlinear 2d case we need to establish an $L^\infty((0,T)\times\Omega)$ estimate of $u$ within the coefficient $c^{-2}-2\gamma u$ in order to guarantee nondegeneracy. We will do so via the embedding 
$H^2((0,T)\times\Omega)\to L^\infty((0,T)\times\Omega)$, the Poincar\'{e} inequality applied to the domain $(0,T)\times\Omega$ with fixed Cauchy data on the boundary part $\{0\}\times\Omega$, as well as the energy estimate resulting from multiplication of the PDE with $-\Delta u$:
\begin{eqnarray}\label{energyid4}
\lefteqn{\frac12\int_0^t  \|\Delta u\|_{L^2(\Omega)}^2 \, ds + \frac{\beta}{2}\|\Delta u(t)\|_{L^2(\Omega)}^2}
\nonumber\\
&=&\frac{\beta}{2}\|\Delta u_0\|_{L^2(\Omega)}^2
+\int_0^t\Bigl(-\Delta u - \beta \Delta u_t\Bigr)(-\Delta u)\, d\Omega\, ds 
- \frac12\int_0^t  \|\Delta u\|_{L^2(\Omega)}^2 \, ds
\nonumber\\
&=&\frac{\beta}{2}\|\Delta u_0\|_{L^2(\Omega)}^2
+\int_0^t\Bigl(\alpha u_{tt}-f u_t-g\Bigr)\Delta u\, d\Omega\, ds 
- \frac12\int_0^t  \|\Delta u\|_{L^2(\Omega)}^2 \, ds
\nonumber\\
&\leq&\frac{\beta}{2}\|\Delta u_0\|_{L^2(\Omega)}^2
+ \int_0^t \|\alpha u_{tt}\|_{L^2(\Omega)}^2\, ds+2\int_0^t \|f u_t\|_{L^2(\Omega)}^2\, ds+2\int_0^t \|g\|_{L^2(\Omega)}^2\, ds
\nonumber\\
\end{eqnarray}

\subsection{Zero and first order ABCs in 1-d; proof of Theorems \ref{theo:zero1d}, \ref{theo:first1d}}\label{subsecwell1}

We prove well-posedness and boundedness of the enrgy $E_1[u]$ as in \eqref{E1} of the following initial boundary value problems
\begin{equation}
\begin{split}
& (c^{-2}-2\gamma u)u_{tt}-u_{xx}-\beta u_{txx}=2\gamma (u_t)^2\quad\text{in } (0,T)\times\Omega,\\
& u(t=0)=u_0\,,\ u_t(t=0)=u_1\quad\text{in }\Omega,\\
& \left.\sqrt{c^{-2}-2\gamma u}\, u_t \pm u^\beta_x\right|_{x=\pm1}=0\quad\text{at } (0,T)\times\{\pm1\},,
\end{split}
\label{Westervelt_zero_1D_IBVP}
\end{equation}
\begin{equation}
\begin{split}
& (c^{-2}-2\gamma u)u_{tt}-u_{xx}-\beta u_{txx}=2\gamma (u_t)^2\quad\text{in } (0,T)\times\Omega,\\
& u(t=0)=u_0\,,\ u_t(t=0)=u_1\quad\text{in }\Omega,\\
& \left. \sqrt{c^{-2}-2\gamma u}\,u_t\pm u^\beta_x
-\frac{\gamma}{2\sqrt{c^{-2}-2\gamma u}}\left(
u_t u \mp\frac{1}{\sqrt{c^{-2}-2\gamma u}} u^\beta_x u\right)
\right|_{x=\pm1}=0\\
&\hspace*{7cm}\text{at } (0,T)\times\{\pm1\},,\\
\end{split}
\label{Westervelt_first_1D_IBVP}
\end{equation}
with $\Omega=(-1,1)$.

To this end we use a fixed point argument for the operator $\mathcal{T}$ mapping $v\in\mathcal{W}$ to a soluion $u$ of 
\begin{equation}\label{West1d_fix}
\begin{split}
& (c^{-2}-2\gamma v)u_{tt}-u_{xx}-\beta u_{txx}=2\gamma v_tu_t\quad\text{in } (0,T)\times\Omega,\\
& u(t=0)=u_0\,,\ u_t(t=0)=u_1\quad\text{in }\Omega,\\
& \left.\frac{2(c^{-2}-2\gamma v)-\zeta \gamma v}{2(c^{-2}-2\gamma v)+\zeta \gamma v}
\sqrt{c^{-2}-2\gamma v}\, u_t 
\pm u^\beta_{x} 
\right|_{x=\pm1}=0\quad
\text{at } (0,T)\times\{\pm1\},,
\end{split}
\end{equation}
where 
$$ \zeta=\left\{\begin{array}{ll} 
0 & \mbox{ in case of }\eqref{Westervelt_zero_1D_IBVP}\\
1 & \mbox{ in case of }\eqref{Westervelt_first_1D_IBVP}
\end{array}\right.$$
and
\begin{equation}\label{defW}
\begin{aligned}
\mathcal{W}=\{&v\in L^2(0,T;L^2(\Omega)) \ : \ v(t=0)=u_0\,, \ v_t(t=0)=u_1\,, \\
&-\underline{m}\leq v(t,x)\leq \bar{m}\,, \ (t,x)\in(0,T)\times\Omega\\
& \|v_t\|_{C(0,T;L^1(\Omega))}\leq \bar{a}
\,,\ \|v_{tt}\|_{L^2(0,T;L^\infty(\Omega))}\leq \bar{b}
\,,\ \|v_{tt}(\pm1)\|_{L^2(0,T)}\leq \bar{c}
\}
\end{aligned}
\end{equation}
for fixed bounds $0<\bar{m} < \frac{1}{2c^2\gamma }$, $0<\underline{m},\bar{a},\bar{b},\bar{c}$ 
(sufficiently small), and we assume that the initial data also satisfy these bounds
$$
-\underline{m}\leq u_0(x)\leq \bar{m}\,, \ x\in\Omega \,, \  
\|u_1\|_{L^1(\Omega)}\leq \bar{a}\,,\ 
\|u_2\|_{L^2(\Omega)}\leq \bar{b}\,,
$$
where 
$$
u_2=\frac{1}{c^{-2}-2\gamma u_0}\left(u_{0xx}+\beta u_{1xx}+2\gamma (u_1)^2\right)\,.
$$

The energy identity \eqref{energyid1} from the previous section with $\alpha=c^{-2}-2\gamma v$, $f=2\gamma v_t$, $g=0$ yields for 
$E_1[u]$ as in \eqref{E1}
\begin{align*}
&E_1[u](t)
+\beta\int_0^t \|u_{ttx}\|_{L^2(\Omega)}^2\, ds\\
&+\int_0^t \Bigl(|\sqrt{\frac{2\alpha-\zeta \gamma v}{2\alpha+\zeta \gamma v}}\alpha^{1/4}u_{tt}|^2(1)
+|\sqrt{\frac{2\alpha-\zeta \gamma v}{2\alpha+\zeta \gamma v}}\alpha^{1/4}u_{tt}|^2(-1)\Bigr)\, ds\\
&= E_1[u](0)
+\gamma\int_0^t\int_\Omega \Bigl( 3v_t (u_{tt})^2 + 2v_{tt} u_t u_{tt} \Bigr)ds \\
&\leq E_1[u](0) +3\gamma \|v_t\|_{C(0,T;L^1(\Omega))} \int_0^t \|u_{tt}(s)\|_{L^\infty(\Omega)}^2\, ds\\
&\quad+\gamma \sqrt{\int_0^t \|v_{tt}(s)\|_{L^\infty(\Omega)}^2\, ds}  
\Bigl(\|u_t\|_{C(0,T;L^1(\Omega))}^2
+\|u_{tt}(s)\|_{L^\infty(\Omega)}^2\, ds\Bigr)\\
&\quad+\gamma \int_0^t \Bigl(\frac{4(1+\zeta)\alpha^2+8\zeta\alpha\gamma v-\zeta^2\gamma^2 v^2}{
\sqrt{\alpha}(2\alpha+\zeta\gamma v)^2} v_t u_t u_{tt}\Bigr)(1) \, ds
\\
&\quad+\gamma \int_0^t \Bigl(\frac{4(1+\zeta)\alpha^2+8\zeta\alpha\gamma v-\zeta^2\gamma^2 v^2}{
\sqrt{\alpha}(2\alpha+\zeta\gamma v)^2} v_t u_t u_{tt}\Bigr)(-1) \, ds\\
&\leq E_1[u](0) 
+\gamma (6\|v_t\|_{C(0,T;L^1(\Omega))}+2\|v_{tt}\|_{L^2(0,T;L^\infty(\Omega))}) 
\Bigl(\int_0^t |u_{tt}(s,-1)|^2 \, ds + 2\int_0^t \|u_{ttx}\|_{L^2(\Omega)}^2\, ds \Bigr)\\
&\quad+\gamma \left\|\frac{2\alpha^2+2\zeta(\alpha+\gamma v)^2}{\sqrt{\alpha}(2\alpha+\zeta\gamma v)^2}\right\|_{C([0,T]\times\overline{\Omega})}
\Bigl\{||v_t(1)||_{C[0,t]}\int_0^t\Bigl(|u_t(s,1)|^2+|u_{tt}(s,1)|^2\Bigr) \, ds\\
&\hspace*{5.5cm}+||v_t(-1)||_{C[0,t]}\int_0^t\Bigl(|u_t(s,-1)|^2+|u_{tt}(s,-1)|^2\Bigr) \, ds\Bigr\}\\
&\leq E_1[u](0) \\
&\quad+ (2T^2(|u_1(1)|+ \sqrt{T} \|v_{tt}(1)\|_{L^2(0,t)}) |u_1(1)|
+ 2T^2(|u_1(-1)|+ \sqrt{T} \|v_{tt}(-1)\|_{L^2(0,t)}) |u_1(-1)|)\\
&\quad\qquad \cdot\left\|\frac{2\alpha^2+2\zeta(\alpha+\gamma v)^2}{\sqrt{\alpha}(2\alpha+\zeta\gamma v)^2}\right\|_{C([0,T]\times\overline{\Omega})}\\
&\quad+2\gamma (3\|v_t\|_{C(0,T;L^1(\Omega))}+\|v_{tt}\|_{C(0,T;L^\infty(\Omega))}) 
\Bigl(\int_0^t |u_{tt}(s,-1)|^2 \, ds + 2\int_0^t \|u_{ttx}\|_{L^2(\Omega)}^2\, ds \Bigr)\\
&\quad+\gamma (1+2T^3)\left\|\frac{2\alpha^2+2\zeta(\alpha+\gamma v)^2}{\sqrt{\alpha}(2\alpha+\zeta\gamma v)^2}\right\|_{C([0,T]\times\overline{\Omega})}
\Bigl\{(|u_1(1)|+ \sqrt{T} \|v_{tt}(1)\|_{L^2(0,t)})\int_0^t|u_{tt}(s,1)|^2 \, ds
\\&\quad\hspace*{6cm}
+(|u_1(-1)|+ \sqrt{T} \|v_{tt}(-1)\|_{L^2(0,t)})\int_0^t|u_{tt}(s,-1)|^2 \, ds
\Bigr\}
\end{align*}
where we have used
\begin{eqnarray}\label{embedding1}
\lefteqn{\int_0^t \|u_{tt}(s)\|_{L^\infty(\Omega)}^2\, ds 
\ \leq \ \int_0^t \|u_{tt}(s)\|_{C(\overline{\Omega})}^2\, ds}\nonumber\\
&=&\int_0^t \sup_{x\in\overline{\Omega}} |u_{tt}(s,x)|^2\, ds 
= \int_0^t \sup_{x\in\overline{\Omega}} |u_{tt}(s,-1)+\int_{-1}^x u_{ttx}(s,\xi)\, d\xi|^2\, ds\nonumber\\
&\leq& 2 \Bigl(\int_0^t |u_{tt}(s,-1)|^2 \, ds + 2\int_0^t \|u_{ttx}\|_{L^2(\Omega)}^2\, ds \Bigr)
\end{eqnarray}
and
\begin{eqnarray}\label{embedding2}
|v_t(t,1)| &=& |u_1(1)+\int_0^t v_{tt}(s,1)\, ds| 
\leq |u_1(1)|+ \sqrt{T} \sqrt{\int_0^t |v_{tt}(s,1)|^2\, ds} \nonumber\\
\int_0^t|u_t(s,1)| \, ds &=& \int_0^t|u_1(1)+\int_0^s u_{tt}(\sigma,1)\,d\sigma\, ds|\nonumber\\
& \leq& T|u_1(1)|+ T^{3/2} \sqrt{\int_0^t |u_{tt}(s,1)|^2\, ds}\,.\nonumber\\
\end{eqnarray}
Since $\overline{\alpha}=c^{-2}+2\gamma\underline{m}\geq\alpha=c^{-2}-2\gamma v\geq c^{-2}-2\gamma\bar{m}=:\underline{\alpha}>0$,  
$\|v_t\|_{C(0,T;L^1(\Omega))}\leq\bar{a}$, $\|v_{tt}\|_{L^2(0,T;L^\infty(\Omega))}\leq \bar{b}$, $\|v_{tt}(\pm1)\|_{L^2(0,T)}\leq \bar{c}$,
and 
\begin{equation}\label{estctil}
\left\|\frac{2\alpha^2+2\zeta(\alpha+\gamma v)^2}{\sqrt{\alpha}(2\alpha+\zeta\gamma v)^2}\right\|_{C([0,T]\times\overline{\Omega})}\leq\frac{2\overline{\alpha}^2+2\zeta(\overline{\alpha}-\gamma\underline{m})^2}{\sqrt{\underline{\alpha}}(2c^{-2}-(4-\zeta)\gamma \overline{m})^2}
=:\tilde{C}(\underline{m},\overline{m},\underline{\alpha},\overline{\alpha})\,,
\end{equation}
this yields
\begin{eqnarray*}
\lefteqn{E_1[u](t)+\tilde{\beta} \int_0^t \|u_{ttx}\|_{L^2(\Omega)}^2\, ds
+\tilde{\underline{\alpha}}\int_0^t \Bigl(|u_{tt}|^2(1)+|u_{tt}|^2(-1)\Bigr)\, ds}\\
&\leq& E_1[u](0)+
2(2T^2(\max\{|u_1(1)|,|u_1(-1)|\}+ \sqrt{T} \bar{c}) \\
&&\qquad \qquad \times\max\{|u_1(1)|,|u_1(-1)|\}
\tilde{C}(\underline{m},\overline{m},\underline{\alpha},\overline{\alpha})
\end{eqnarray*}
with 
$\tilde{\beta}=\beta-12\gamma\bar{a}-4\gamma\bar{b}>0$, \\
$\tilde{\underline{\alpha}}=\sqrt{\underline{\alpha}}- 6\gamma\bar{a}- 2\gamma\bar{b}
\gamma (1+2T^3) \tilde{C}(\underline{m},\overline{m},\underline{\alpha},\overline{\alpha})
(\max\{|u_1(1)|,|u_1(-1)|\}+\sqrt{T}\bar{c})>0$ \\
for $\bar{a}, \bar{c}$ sufficiently small. 
The definition of $E_1[u]$ as well as the $C([0,T]\times\overline{\Omega})$ estimate
\begin{eqnarray}\label{Linfty1d}
|u(t,x)|&=& |u_0(x)+\int_0^t u_t(s,x)\, ds| = |u_0(x)+\int_0^t (u_1(x)+\int_0^s u_{tt}(\sigma,x)\, d\sigma|\nonumber\\
&\leq& \|u_0\|_{L^\infty(\Omega)} + T \|u_1\|_{L^\infty(\Omega)} +T^{3/2}\sqrt{\int_0^t |u_{tt}(s,x)|^2\, ds} 
\end{eqnarray}
together with \eqref{embedding1} after possibly decreasing $\bar{a}$ allows us to conclude \\
$-\underline{m}\leq u(t,x)\leq \bar{m}\,, \ (t,x)\in(0,T)\times\Omega$,\\
$\|u_t\|_{C(0,T;L^1(\Omega))}\leq \bar{a}$, 
$\|u_{tt}\|_{L^2(0,T;L^\infty(\Omega))}\leq \bar{b}$,
$\|u_{tt}(\pm1)\|_{L^2(0,T)}\leq \bar{c}$,\\
i.e., altogether we have $u\in \mathcal{W}$.
Note that the appearance of constants depending on $T$ in \eqref{embedding2}, \eqref{Linfty1d} prevents us from showing global existence and exponential decay here.

To show that $\mathcal{T}$ is a contraction on $\mathcal{W}$, we use the fact that for $v^1,v^2\in\mathcal{W}$, and  $u^i = \mathcal{T} v_i$, $i=1,2$, the function $\hat{u}=u^1 -u^2$ solves the following problem ($\hat{v} = v^1 - v^2$)
\begin{equation}\label{contr}
\begin{split}
&(c^{-2}-2\gamma v^1) \hat{u}_{tt} -\hat{u}_{xx} - \beta \hat{u}_{txx} = 
2\gamma v^1_t \hat{u}_t+ 2\gamma \hat{v}_t u^2_t + 2\gamma \hat{v} u^2_{tt} 
\quad\text{in } (0,T)\times\Omega,\\
& \hat{u}(t=0)=0\,,\ \hat{u}_t(t=0)=0\quad\text{in }\Omega,\\
& \left.\frac{2(c^{-2}-2\gamma v^1)-\zeta \gamma v^1}{2(c^{-2}-2\gamma v^1)+\zeta \gamma v^1}
\sqrt{c^{-2}-2\gamma v^1}\,\hat{u}_t\pm(\hat{u}_x+\beta \hat{u}_{tx})
+\tilde{c}^{-1} \hat{v}\right|_{x=0}=0\\
&\hspace*{7cm}\quad\text{at } (0,T)\times\{\pm1\}\,,
\end{split}
\end{equation}
where
\begin{eqnarray*}
\tilde{c}^{-1}(t,x)
&=& \frac{\frac{2(c^{-2}-2\gamma v^2)-\zeta \gamma v^2}{2(c^{-2}-2\gamma v^2)+\zeta \gamma v^2}
\sqrt{c^{-2}-2\gamma v^2} 
-\frac{2(c^{-2}-2\gamma v^1)-\zeta \gamma v^1}{2(c^{-2}-2\gamma v^1)+\zeta \gamma v^1}
\sqrt{c^{-2}-2\gamma v^1}
}{v^2(t,x)-v^1(t,x)} u^2_t(t,x)\\
&=&-\gamma\int_0^1 \frac{4(1+\zeta)(\alpha^\theta)^2+8\zeta\alpha^\theta\gamma v^\theta-\zeta^2\gamma^2 
(v^\theta)^2}{
\sqrt{\alpha}(2\alpha^\theta+\zeta\gamma v^\theta)^2}
d\theta \, u^2_t(t,x)\,,
\end{eqnarray*}
where $v^\theta=v^2(t,x)+\theta \hat{v}(t,x))$, $\alpha^\theta=c^2-2\gamma v^\theta$.
Hence using \eqref{energyid0} with $\alpha=c^{-2}-2\gamma v^1$, $f=2\gamma v^1_t $, $g=2\gamma \hat{v}_t u^2_t + 2\gamma \hat{v} u^2_{tt}$ we obtain for 
$$
E_0[\hat{u}](t)=
\frac12 \Bigl(\|\sqrt{\alpha(t)} \hat{u}_t(t)\|_{L^2(\Omega)}^2 + \|\nabla \hat{u}(t)\|_{L^2(\Omega)}^2 \Bigr)
$$
the estimate
\begin{eqnarray*}
\lefteqn{E_0[\hat{u}](t)
+\beta\int_0^t \|\hat{u}_{tx}\|_{L^2(\Omega)}^2\, ds}\\
&&+\int_0^t \Bigl(
|\sqrt{\frac{2\alpha-\zeta \gamma v^1}{2\alpha+\zeta \gamma v^1}}\alpha^{1/4}\hat{u}_t(s,1)|^2
+|\sqrt{\frac{2\alpha-\zeta \gamma v^1}{2\alpha+\zeta \gamma v^1}}\alpha^{1/4}\hat{u}_t(s,-1)|^2\Bigr)\, ds\\
&\leq& E_0[\hat{u}](0) + 2\gamma\int_0^t\int_\Omega v^1_t (\hat{u}_t)^2 \, dx\, ds\\
&&+ 2\gamma \int_0^t\int_\Omega (\hat{v}_t u^2_t + \hat{v} u^2_{tt}) \hat{u}_t\, dx\, ds
-\int_0^t (\tilde{c}^{-1} \hat{v} \hat{u}_t)(s,-1)+\tilde{c}^{-1} \hat{v} \hat{u}_t)(s,1))\, ds  \\
&\leq& 2\gamma \|v^1_t\|_{C(0,T;L^1(\Omega))} \int_0^t \|\hat{u}_t(s)\|_{L^\infty(\Omega)}^2\, ds \\
&&+ 2\gamma \|u^2_t\|_{C(0,T;L^1(\Omega))}\int_0^t \|\hat{u}_t(s)\|_{L^\infty(\Omega)}\|\hat{v}_t(s)\|_{L^\infty(\Omega)}\, ds
\\
&&+ 2\gamma \|\hat{v}\|_{C([0,T]\times\overline{\Omega})} \int_0^t \|u^2_{tt}(s)\|_{L^1(\Omega)} \|\hat{u}_t(s)\|_{L^\infty(\Omega)}\, ds\\
&&+ \|\hat{v}\|_{C([0,T]\times\overline{\Omega})}
\int_0^t \Bigl(|\tilde{c}^{-1}(s,1)|\,|\hat{u}_t(s,1)|+|\tilde{c}^{-1}(s,-1)|\,|\hat{u}_t(s,-1)|\Bigr)\, ds\\
&\leq& 2\gamma \bar{a} \int_0^t \|\hat{u}_t(s)\|_{L^\infty(\Omega)}^2\, ds \\
&&+ \gamma\bar{a} 
\Bigl(\int_0^t \|\hat{u}_t(s)\|_{L^\infty(\Omega)}^2+\int_0^t\|\hat{v}_t(s)\|_{L^\infty(\Omega)}^2\, ds\Bigr)\\
&&+ \gamma \bar{b}
\Bigl(\|\hat{v}\|_{C([0,T]\times\overline{\Omega})}^2+\int_0^t\|\hat{u}_t(s)\|_{L^\infty(\Omega)}^2\, ds\Bigr)\\
&&+ \gamma \bar{c} 
\frac{2\overline{\alpha}^2+2\zeta(\overline{\alpha}-\gamma\underline{m})^2}{\sqrt{\underline{\alpha}}(2c^{-2}-(4-\zeta)\gamma \overline{m})^2}
\Bigl(\frac12\|\hat{v}\|_{C([0,T]\times\overline{\Omega})}^2
+\int_0^t \Bigl(|\hat{u}_t(s,1)|^2+|\hat{u}_t(s,-1)|^2\Bigr)\, ds\Bigr)\,,
\end{eqnarray*}
where we have used $E_0[\hat{u}](0)=0$ as well as an estimate similar to \eqref{estctil} together with 
$\|u^2_t(\pm1)\|_{L^2(0,T)}\leq \bar{c}$ to estimate $\|\tilde{c}^{-1}(\pm1)\|_{L^2(0,T)}$.
Since $\alpha^{1/4}\geq \sqrt[4]{c^{-2}-2\gamma\bar{m}}$, and by estimate \eqref{embedding1} for 
$\hat{u}_t, \hat{v}_t$ in place of $u_{tt}$, we arrive at an estimate of the form 
\begin{eqnarray*}
\lefteqn{\max\{\|E_0[\hat{u}]\|_{C[0,T]}\,,\ \|\hat{u}(\pm1)\|_{L^2(0,T)}\} }\nonumber\\
\leq \ &&
C(\underline{m}, \bar{m}, \bar{a}, \bar{b}, \bar{c}) \Bigl(
\max\{\|E_0[\hat{u}]\|_{C[0,T]}\,,\ \|\hat{u}(\pm1)\|_{L^2(0,T)}\,, \\
&& \hspace*{5cm}
\|E_0[\hat{v}]\|_{C[0,T]}\,,\ \|\hat{v}(\pm1)\|_{L^2(0,T)}\}\Bigr)
\end{eqnarray*}
with a constant $C(\underline{m}, \bar{m}, \bar{a}, \bar{b}, \bar{c})$ that can be made small by
$\bar{a},\bar{b}$ sufficiently small. Hence, we achieve contractivity of  
$\mathcal{T}$ on $\mathcal{W}$ with respect to the norm induced by $\max\{\|E_0[\hat{u}]\|_{C[0,T]}\,,\ \|\hat{u}(\pm1)\|_{L^2(0,T)}\}$.

Thus, using Banach's Contraction Principle, we have shown Theorems \ref{theo:zero1d}, \ref{theo:first1d}.

\begin{remark}\label{rem:timedepcoeff}
It is readily checked that replacing $\sqrt{c^{-2}-2\gamma u}$ by a constant $c^{-1}$, we would get rid of a couple of higher derivative terms on the boundary and 
end up with energy estimates enabling even global in time wellposedness even with $\beta=0$, cf. \cite{K_amo}. Thus, for obtaining enhanced approximation by taking 
into account the full time and space dependence of the coefficient $\sqrt{c^{-2}-2\gamma u}$ in \eqref{Westervelt_zero_1D_IBVP} we pay the price of losing global 
in time well-posedness and needing strong damping. 
\end{remark}

\begin{remark}\label{rem:maxparreg}
For $\beta>0$ it is possible to make use of maximal parabolic regularity to prove more general existence results in $L^p$ spaces even under less restrictive 
assumptions on the regularity of the inital data in case of pure Dirichlet boundary conditions cf. \cite{Wilke}. However, it seems to be at least not really 
straightforward to carry over these techniques to absorbing boundary conditions.  
\end{remark}

\subsection{Zero and first order ABCs in 2-d; proof of Theorems \ref{theo:zero2d}, \ref{theo:first2d}}\label{subsecwell2}

Again we use Banach's Contraction Principle, where the self mapping property on balls $\mathcal{W}$ with respect to the combined energy 
$$
E_3[v](t)=
\frac12 \Bigl(\|\sqrt{\alpha(t)} u_{tt}(t)\|_{L^2(\Omega)}^2 + \| \nabla u_t(t)\|_{L^2(\Omega)}^2 
+\lambda\beta\|\Delta u(t)\|_{L^2(\Omega)}^2
\Bigr)
$$
in case of zero order ABC 
and the combined energy 
\begin{align*}
E_4[u](t)=&\Bigl(\|\alpha(t) u_{tt}(t)\|_{L^2(\Omega)}^2 + \|\sqrt{\alpha(t)} \nabla u_t(t)\|_{L^2(\Omega)}^2\\ 
&+\|\sqrt{\alpha(t)} u^\beta_{tx}(t)\|_{L^2(\Omega)}^2 + \|\nabla u^\beta_x(t)\|_{L^2(\Omega)}^2
+\lambda\beta\|\Delta u(t)\|_{L^2(\Omega)}^2
\Bigr)
\end{align*}
in case of first order ABC, with an appropriately chosen factor $\lambda>0$. Moreover, to include pointwise bounds $-\underline{m}, \bar{m}$ on $u$ into the definition of $\mathcal{W}$ for avoiding degeneracy, we make use of the fact that in both cases obviously an estimate of the form
$$ 
\|u\|_{L^\infty((0,T)\times\Omega)}\leq C \sup_{t\in(0,T)} E_j(t) \quad j\in\{3,4\}
$$
with $C$ possibly depending on $T$ holds.
Accordingly, we will use combination of the energy identities \eqref{energyid1}, \eqref{energyid4} and \eqref{energyid2}, \eqref{energyid3}, \eqref{energyid4}, respectively for showing the self mapping property in case of zero and first order ABCs, respectively.
Contractivity in case of zero order ABC will rely on the lower order energy identity \eqref{energyid0}, applied to the initial boundary value problem that holds for the difference between two solutions of the linearized problem. In case of first order ABCs we will have to use a higher order energy identity also for contractivity in order to take into account the tangential derivative terms. This is the only part of the proof that we will provide explicitely here, since the rest (self-mapping for zero and first order ABCs, contraction for zero order ABCs) goes very much along the lines of the proofs in \cite{Clason2009}, \cite{KL09}:

To show  that the operator $\mathcal{T}$ mapping $v\in \mathcal{W}$ to a solution $u$ of 
\begin{equation}
\begin{split}
& (c^{-2}-2\gamma v)u_{tt}-u_{xx}-u_{yy}-\beta u_{txx}-\beta u_{tyy}=2\gamma v_tu_t\quad\text{in } (0,T)\times\Omega,\\
& u(t=0)=u_0\,,\ u_t(t=0)=u_1\quad\text{in }\Omega,\\
& \left( 
\sqrt{{\alpha^v}}\frac{u_{tt}+u^\beta_{tt}}{2}+u^\beta_{tn}
-\frac{1}{2\sqrt{{\alpha^v}}} (u^\beta_{\vartheta\vartheta}+\beta u^\beta_{\vartheta\vartheta})
\right.\\
&\left.
\left.-\frac{\gamma}{2\sqrt{\alpha}}
\left(u_t -\frac{1}{\sqrt{{\alpha^v}}} u_n \right)v_t
-\frac{\gamma}{2(\alpha^v)^{3/2}}
\left(\frac12 u_t -\frac{1}{\sqrt{{\alpha^v}}} u_n \right)\int_0^\cdot (v^\beta_{\vartheta\vartheta}+\beta v^\beta_{\vartheta\vartheta})\, dt\right)
\right|_{\partial\Omega}=0\\
&\hspace*{7cm}\text{at } (0,T)\times\partial\Omega\,,\\
\end{split}
\label{Westervelt_first_2D_IBVP_intro}
\end{equation}
with $\alpha^v=c^{-2}-2\gamma v$ is a contraction on 
\begin{equation}\label{defW_first2d} 
\begin{aligned}
\mathcal{W}=\{&v\in L^2(0,T;L^2(\Omega)) \ : \ v(t=0)=u_0\,, \ v_t(t=0)=u_1\,, \\
&-\underline{m}\leq v(t,x)\leq \bar{m}\,, \ (t,x)\in(0,T)\times\Omega\\
& \|E_4[v]\|_{C[0,T]}\leq\bar{a}^2\,,\nonumber\\
& \|\sqrt{\alpha^v} \nabla v_{tt}\|_{L^2(0,T;L^2(\Omega))}\,, \ 
\|\nabla v^\beta_{tx}\|_{L^2(0,T;L^2(\Omega))} \,, \ 
\|\Delta v\|_{L^2(0,T;L^2(\Omega))} \ \leq \bar{b}\,,\nonumber\\
& \|(\alpha^v)^{1/4} v^\beta_{tn}\|_{L^2(0,T;L^2(\partial\Omega))}\,, \
\|\alpha^v (v+v^\beta)_{tt}  - (v^\beta_{\vartheta\vartheta}+\beta v^\beta_{t\vartheta\vartheta})\|_{L^2(0,T;L^2(\partial\Omega))} \ \leq \bar{c}
\}\,,
\end{aligned}
\end{equation}
we use the fact that for $v^1,v^2\in\mathcal{W}$, and  $u^i = \mathcal{T} v_i$, $i=1,2$, the function $\hat{u}=u^1 -u^2$ solves the following problem ($\hat{v} = v^1 - v^2$)
\begin{equation}\label{contr}
\begin{split}
&\alpha^{v^1} \hat{u}_{tt} -\hat{u}_{xx} - \beta \hat{u}_{txx} = 
2\gamma v^1_t \hat{u}_t+ 2\gamma \hat{v}_t u^2_t + 2\gamma \hat{v} u^2_{tt} 
\quad\text{in } (0,T)\times\Omega,\\
& \hat{u}(t=0)=0\,,\ \hat{u}_t(t=0)=0\quad\text{in }\Omega,\\
&  \sqrt{{\alpha^{v^1}}}\frac{\hat{u}_{tt}+\hat{u}^\beta_{tt}}{2}+\hat{u}^\beta_{tn}
-\frac{1}{2\sqrt{{\alpha^{v^1}}}} (\hat{u}^\beta_{\vartheta\vartheta}+\beta \hat{u}^\beta_{\vartheta\vartheta})
\\
& -\frac{\gamma}{2\sqrt{\alpha^{v^1}}}
\left(\hat{u}_t -\frac{1}{\sqrt{{\alpha^{v^1}}}} \hat{u}_n \right){v^1}_t
-\frac{\gamma}{2(\alpha^{v^1})^{3/2}}
\left(\frac12 \hat{u}_t -\frac{1}{\sqrt{{\alpha^{v^1}}}} \hat{u}_n \right)\int_0^\cdot ({v^1}^\beta_{\vartheta\vartheta}+\beta {v^1}^\beta_{\vartheta\vartheta})\, dt\\
& +\left. \tilde{c}^{-1} \hat{v}\right|_{x=0}=0
\quad\text{at } (0,T)\times\{\pm1\},,
\end{split}
\end{equation}
where
\begin{eqnarray*}
\tilde{c}^{-1}
&=& \frac{\sqrt{\alpha^{v^2}}-\sqrt{\alpha^{v^1}}}{v^2-v^1} \frac{u^2_{tt}+{u^2}^\beta_{tt}}{2}\\
&&-\frac{\frac{1}{2\sqrt{{\alpha^{v^2}}}}-\frac{1}{2\sqrt{{\alpha^{v^1}}}}}{v^2-v^1}
 ({u^2}^\beta_{\vartheta\vartheta}+\beta {u^2}^\beta_{\vartheta\vartheta})\\
&&-\frac{\frac{\gamma}{2\sqrt{\alpha^{v^2}}}+\frac{\gamma}{4(\alpha^{v^2})^{3/2}}
-\frac{\gamma}{2\sqrt{\alpha^{v^1}}}+\frac{\gamma}{4(\alpha^{v^1})^{3/2}}}{v^2-v^1} u^2_t\\
&&+\frac{\frac{\gamma}{2\alpha^{v^2}}+\frac{\gamma}{2(\alpha^{v^2})^2}
\int_0^\cdot ({v^2}^\beta_{\vartheta\vartheta}+\beta {v^2}^\beta_{\vartheta\vartheta})\, dt
-\frac{\gamma}{2\alpha^{v^1}}+\frac{\gamma}{2(\alpha^{v^1})^2}
\int_0^\cdot ({v^1}^\beta_{\vartheta\vartheta}+\beta {v^1}^\beta_{\vartheta\vartheta})\, dt
}{v^2-v^1} u^2_n\,.
\end{eqnarray*}
Hence using \eqref{energyid23} with $\alpha=\alpha^{v^1}=c^{-2}-2\gamma v^1$, $f=2\gamma v^1_t $, $g=2\gamma \hat{v}_t u^2_t + 2\gamma \hat{v} u^2_{tt}$ we obtain for the energy
\begin{align*}
E_2[\hat{u}](t)= 
\frac12 \Bigl(&\|\alpha(t) \hat{u}_{tt}(t)\|_{L^2(\Omega)}^2+\|\sqrt{\alpha(t)} \nabla \hat{u}_t(t)\|_{L^2(\Omega)}^2\\
&+ \|\sqrt{\alpha(t)} \hat{u}^\beta_{tx}(t)\|_{L^2(\Omega)}^2 + \|\nabla \hat{u}^\beta_x(t)\|_{L^2(\Omega)}^2 \Bigr)
\end{align*}
and the interior dissipation
$$
D_2[\hat{u}](t)= 
\beta\Bigl(\|\sqrt{\alpha} \nabla \hat{u}_{tt}\|_{L^2(\Omega)}^2+\|\nabla \hat{u}^\beta_{tx}\|_{L^2(\Omega)}^2\Bigr)
$$
\begin{eqnarray}\label{estE2}
\lefteqn{E_2[\hat{u}](t)
+\int_0^t D_2[\hat{u}](s)\, ds +2\int_0^t\|\sqrt{\alpha} \hat{u}^\beta_{tn}\|_{L^2(\partial\Omega)} \, ds
= \ E_2[\hat{u}](0) }\nonumber\\
&&+\int_0^t \int_\Omega \Bigl( 
\alpha f (\hat{u}_{tt})^2 
+\frac12\alpha_t |\nabla \hat{u}_t|^2 
+(f-\tfrac{\alpha_t}{2}) (\hat{u}^\beta_{tx})^2
+(\beta f_t +\alpha_t)\hat{u}_{tx}\hat{u}^\beta_{tx}
\nonumber\\
&& \qquad\qquad
-\beta\alpha_{tx} \hat{u}_{tt}\hat{u}^\beta_{tx}
-\alpha_x \hat{u}^\beta_{tt}\hat{u}^\beta_{tx}
+\alpha f_t \hat{u}_t \hat{u}_{tt} 
- \hat{u}_{tt} \nabla \alpha \nabla \hat{u}^\beta_t 
\nonumber\\
&& \qquad\qquad
+f_x \hat{u}^\beta_t \hat{u}^\beta_{tx} 
+\beta f_{tx} \hat{u}_t \hat{u}^\beta_{tx}
+ \alpha g_t \hat{u}_{tt}+(g_x+\beta g_{tx})\hat{u}^\beta_{tx}
\Bigr) \, d\Omega\, ds
\nonumber\\
&&+\int_0^t\int_{\partial\Omega} \Bigl(
\gamma
\left(\hat{u}_t -\frac{1}{\sqrt{{\alpha}}} \hat{u}_n \right){v^1}_t
+\frac{\gamma}{\alpha}
\left(\frac12 \hat{u}_t -\frac{1}{\sqrt{{\alpha}}} \hat{u}_n \right)\int_0^\cdot ({v^1}^\beta_{\vartheta\vartheta}+\beta {v^1}^\beta_{\vartheta\vartheta})\, dt
-\tilde{c}^{-1} \hat{v}
\nonumber\\
&&\qquad\qquad-f\hat{u}^\beta_t-\beta f_t\hat{u}_t
-g-\beta g_t
+\beta\alpha_t \hat{u}_{tt}\Bigr) \hat{u}^\beta_{tn} \, d\Gamma\, ds
\nonumber\\
&=& 2\gamma \int_0^t \int_\Omega \Bigl\{ 
\alpha v^1_t (\hat{u}_{tt})^2 
-\tfrac12 v^1_t |\nabla \hat{u}_t|^2 
+(v^1_t+\tfrac12 v^1_t) (\hat{u}^\beta_{tx})^2
\nonumber\\
&& \qquad\qquad
+(\beta v^1_{tt} -v^1_t)\hat{u}_{tx}\hat{u}^\beta_{tx}
+\beta v^1_{tx} \hat{u}_{tt}\hat{u}^\beta_{tx}
+v^1_x \hat{u}^\beta_{tt}\hat{u}^\beta_{tx}
\nonumber\\&& \qquad\qquad
+\alpha v^1_{tt} \hat{u}_t \hat{u}_{tt} 
+ \hat{u}_{tt} \nabla v^1 \nabla \hat{u}^\beta_t 
+v^1_{tx} \hat{u}^\beta_t \hat{u}^\beta_{tx} 
+\beta v^1_{ttx} \hat{u}_t \hat{u}^\beta_{tx}
\nonumber\\&&\qquad\qquad 
+ \alpha (\hat{v}_t u^2_t + \hat{v} u^2_{tt})_t \hat{u}_{tt}
+((\hat{v}_t u^2_t + \hat{v} u^2_{tt})_x+\beta (\hat{v}_t u^2_t + \hat{v} u^2_{tt})_{tx})\hat{u}^\beta_{tx}
\Bigr\} \, d\Omega\, ds
\nonumber\\
&&+ \gamma \int_0^t\int_{\partial\Omega} \Bigl[
\left(\hat{u}_t -\tfrac{1}{\sqrt{{\alpha}}} \hat{u}_n \right){v^1}_t
+\frac{1}{\alpha}
\left(\tfrac12 \hat{u}_t -\tfrac{1}{\sqrt{{\alpha}}} \hat{u}_n \right)\int_0^\cdot ({v^1}^\beta_{\vartheta\vartheta}+\beta {v^1}^\beta_{\vartheta\vartheta})\, dt
-\tilde{c}^{-1} \hat{v}
\nonumber\\&&\qquad\qquad
-2v^1_t\hat{u}^\beta_t-2\beta v^1_{tt}\hat{u}_t
-2(\hat{v}_t u^2_t + \hat{v} u^2_{tt})-\beta 2(\hat{v}_t u^2_t + \hat{v} u^2_{tt})_t
-2\beta v^1_t \hat{u}_{tt}\Bigr] \hat{u}^\beta_{tn} \, d\Gamma\, ds
 \,,
\nonumber\\
&& \qquad\qquad
\end{eqnarray}
where we have used $E_2[\hat{u}](0)=0$.
Let us first consider the terms within the curly braces $\{\ldots\}$ under the integral over $\Omega$ and $(0,t)$ on the right hand side of \eqref{estE2}. After bounding the $L^\infty$ norm of the $\alpha$ factors by $\overline{\alpha}=c^{-2}+2\gamma\underline{m}$ we see that they are all of the form 
$$ q\cdot\phi\cdot\psi \ \mbox{ with } \ 
q\in \{\d v^1,\d u^2\}
\quad 
\phi,\psi \in \{\partial \hat{u},\partial \hat{v}\}
\,, 
$$
where $\partial$ is a combination of differential operators $(\mbox{id}+\beta\partial_t), \ \partial_x, \ \partial_t, \ \nabla$.
Hence the time and space integrals of these terms can be estimated by products of the form
\begin{eqnarray*}
&&\|f^1\|_{C(0,T;L^2(\Omega)}\cdot\|f^2\|_{L^2(0,T;L^4(\Omega)}\cdot\|f^3\|_{L^2(0,T;L^4(\Omega)}\mbox{ or }\\
&&\|f^1\|_{L^2(0,T;L^2(\Omega)}\cdot\|f^2\|_{C(0,T;L^4(\Omega)}\cdot\|f^3\|_{L^2(0,T;L^4(\Omega)}\\
&& \mbox{where } f^1,f^2,f^3 \mbox{ is an appropriate permutation of }q,\phi,\psi\,.
\end{eqnarray*}
It is readily checked that since $v^1,u^2\in\mathcal{W}$ as defined in \eqref{defW_first2d}, all $q$ factors can be bounded by constants depending on 
$\underline{m}, \bar{m}, \bar{a}, \bar{b}, \bar{c}$ that can be made small for small $\bar{a}, \bar{b}, \bar{c}$, and all $\phi$ and $\psi$ factors can be bounded by either the energy norms 
$\|E_2[\hat{v}]\|_{C[0,T]}$, $\|E_2[\hat{u}]\|_{C[0,T]}$ or the  
interior dissipation norms $\|D_2[\hat{v}]\|_{L^2(0,T)}, \|D_2[\hat{u}]\|_{L^2(0,T)}$.
For the terms within the brackets $[\ldots]$ under the integral over $\partial\Omega$ and $(0,t)$ on the right hand side of \eqref{estE2} we have that their squared $L^2(0,t;L^2(\partial\Omega))$ norm is bounded by some constant depending on  
$\underline{m}, \bar{m}, \bar{a}, \bar{b}, \bar{c}$ (which can be made small for small $\bar{a}, \bar{b}, \bar{c}$), multiplied with $\|E_2[\hat{v}]\|_{C[0,T]}$, $\|E_2[\hat{u}]\|_{C[0,T]}$ or $\|D_2[\hat{v}]\|_{L^2(0,T)}, \|D_2[\hat{u}]\|_{L^2(0,T)}$.
Altogether we arrive at an estimate of the form 
\begin{eqnarray*}
\lefteqn{\max\{\|E_2[\hat{u}]\|_{C[0,T]}\,,\ \|D_2[\hat{u}]\|_{L^2(0,T)}\} }\nonumber\\
\leq \ &&
C(\underline{m}, \bar{m}, \bar{a}, \bar{b}, \bar{c}) \Bigl(
\max\{\|E_2[\hat{u}]\|_{C[0,T]}\,,\ \|D_2[\hat{u}]\|_{L^2(0,T)}\,, 
\nonumber\\
&& \qquad\qquad\qquad\qquad
\|E_2[\hat{v}]\|_{C[0,T]}\,,\ \|D_2[\hat{v}]\|_{L^2(0,T)}\}\Bigr)
\end{eqnarray*}
which for small $\bar{a}, \bar{b}, \bar{c}$ gives the desired contractivity estimate.

\section{Numerical results}
\label{sec:numres}
In this section we study the performance of the proposed boundary conditions and compare them with 
the first and second order Engquist--Majda ABCs \cite{EngquistMajda1979} for different setups. 
In what follows, we focus on a horizontal waveguide in one and two dimensions, and study how the accuracy of ABCs is influenced by the angle of incidence in the 2-d case.
Then, we consider the high-intensity focused ultrasound (HIFU) problem with the physical 
parameters typical for simulations of thermotherapy for
human liver cancer and analyze how intensively the solution is contaminated by the reflected waves.
We name the ABCs as $\rm ABC^{d,o}_n$, where the superscripts {\rm d} and {\rm o} indicate 
the space dimension and the order of ABC, while the subscript $n$ takes the value {\rm PS} or {\rm EM} standing for 
the new nonlinear ABC obtained with the pseudo-differential calculus or the Engquist--Majda ABC, respectively.
To approximate system \eqref{Westervelt}-\eqref{BC} in time the standard Newmark scheme is applied~\cite{Hughes1987}.
For space discretization, the finite element method is used.

In order to compare different ABCs, a reference solution $u^*$ is computed in the domain $\Omega^\prime\Supset\Omega$ which is large enough to prevent
the solution in the restricted domain $\Omega$ from being polluted by reflected waves. The studied ABCs are compared in terms of the $l^2$-norm relative
error $\delta=\|u^*-u\|_2/\|u^*\|_2$. 
In all numerical experiments the number of
finite elements per wavelength is set to be 50, and the time step is chosen in such a way as to have $20$ time samples per time period 
for each of the frequencies $\omega=\{100\, {\rm kHz}, 1\, {\rm MHz}\}$. To induce a wave in the domain, we use a monofrequency transducer of the form
$u_n=\sin(2\pi\omega t)$.
The time $t$ as well as the acoustic pressure are normalized to their maximum values.
The physical parameters used in all the computations correspond to those of human liver~\cite{HallajClevelandHynynen2001,ConnorHynynen2002}: $c=1596\ {\rm m\cdot s^{-1}}$, 
$\rho=1050 \ {\rm kg\cdot m^{-3}}$, $B/A=6.8$, $b=2\alpha c^3/(2\pi\omega)^2$, with the acoustic absorption coefficient $\alpha=4.5\ {\rm Np\cdot m^{-1}\cdot MHz^{-1}}$.

\subsection{ABC in 1-d}
\begin{figure}[htp]
  \hspace*{1.3cm}
  \begin{tabular}{cc}
  \hspace*{0.5cm}\large$\delta$ &
  \begin{minipage}{0.9\textwidth}
    \includegraphics[width=3in]{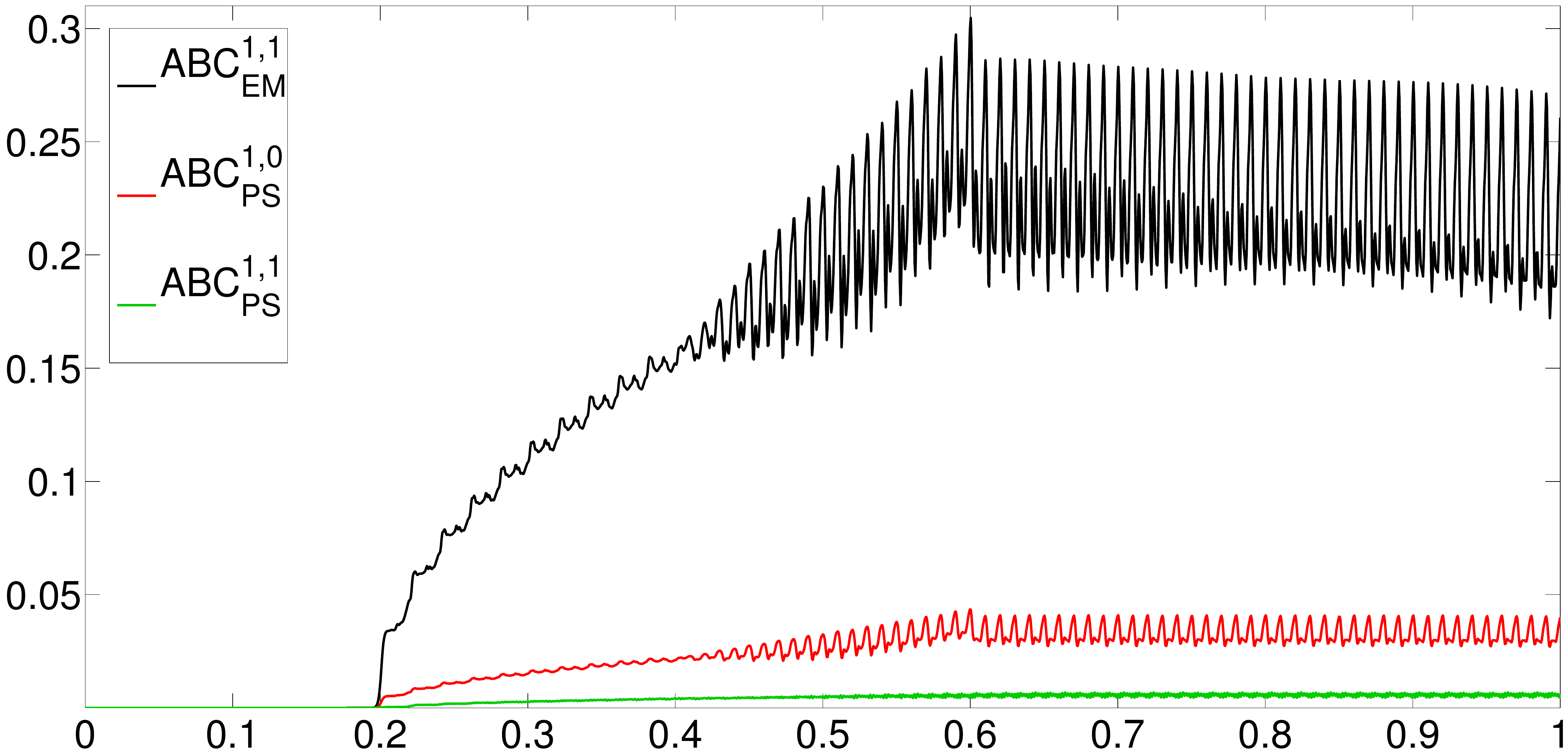}\\
  \end{minipage}\\ 
  \vspace*{-0.7cm}&\\
  & \hspace{-3.4cm}\large $t$  \\
  \end{tabular}  
  \caption{Relative error $\delta$ versus time $t$ for $\rm ABC^{1,o}_n$, $\rm o=\{0,1\}$, $\rm n=\{EM,PS\}$.}
  \label{fig:ABC_1d}
\end{figure}

In this part we compare the zero and first order ABCs on a line segment $x\in[0,3\, {\rm cm}]$.
The transducer with a $100\, {\rm kHz}$ excitation frequency is set at $x=0$ while the ABC is prescribed at $x=3\, {\rm cm}$.
The results of the comparison are presented in figure~\ref{fig:ABC_1d}.

As it can be seen from figure~\ref{fig:ABC_1d}, the behaviour of the boundary conditions brings no surprise: the higher the order of ABC the more accurate solution we have. 
In contrast to $\rm ABC^{1,0}_{PS}$ and $\rm ABC^{1,1}_{PS}$, $\rm ABC^{1,1}_{EM}$ is of much less accuracy. This result is expectable and 
reconfirms the attention one has to pay to the ABCs for the Westervelt equation.

\subsection{ABC in 2-d}
\begin{figure}[htp]
  \hspace*{1.3cm}
  \begin{tabular}{cc}
  \hspace*{0.5cm}\large$\delta$ &
  \begin{minipage}{0.9\textwidth}
    \includegraphics[width=3in]{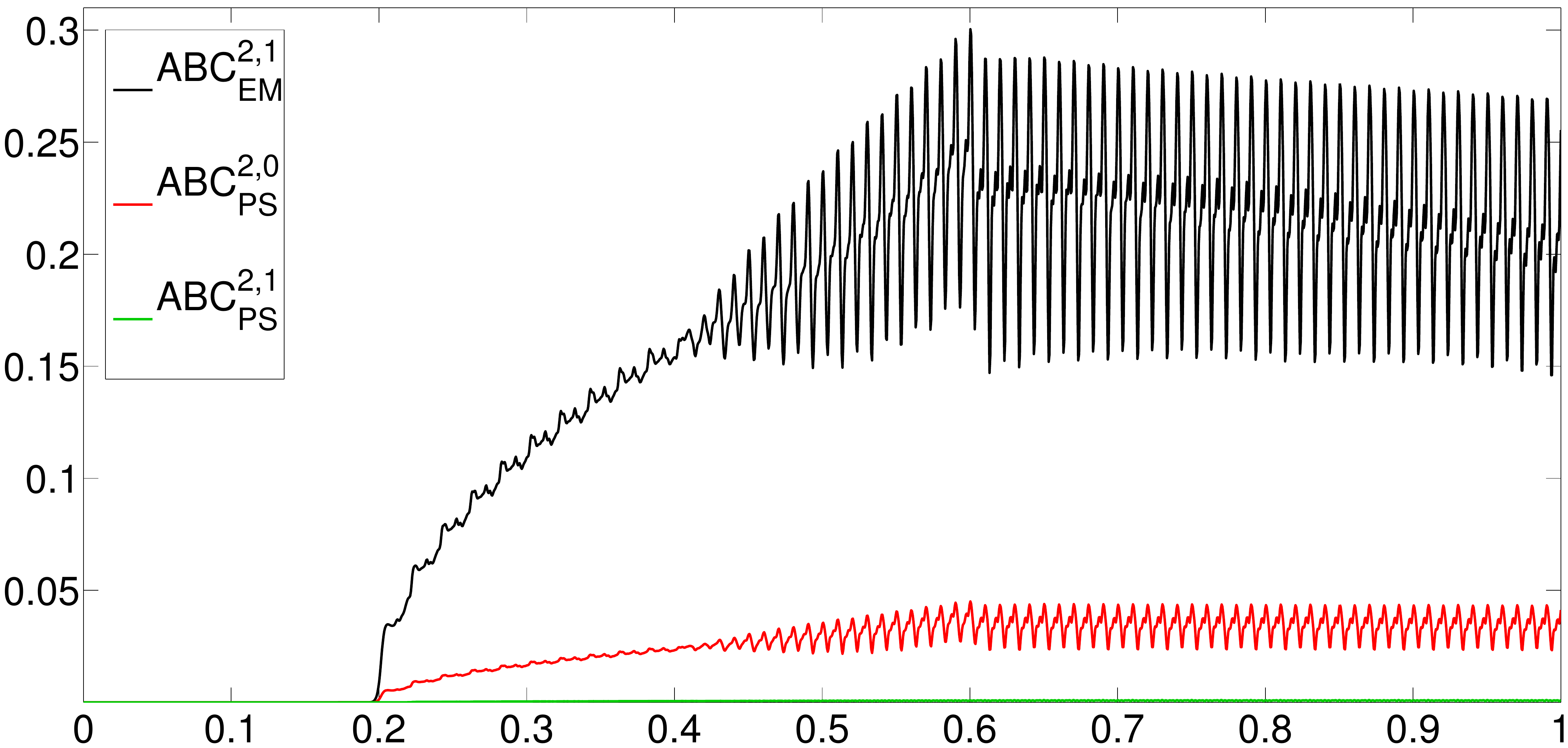}\\
  \end{minipage}\\ 
  \vspace*{-0.7cm}&\\
  & \hspace{-3.4cm}\large $t$  \\
  \end{tabular}  
  \caption{Relative error $\delta$ versus time $t$ for $\rm ABC^{2,o}_n$, $\rm o=\{0,1,2\}$, $\rm n=\{EM,PS\}$.}
  \label{fig:ABC_2d_1}
\end{figure}

The zero order versions of the proposed ABCs, as well as the Engquist--Majda ABCs, provide the best absorption of the wave hitting the boundary at normal incidence. The higher order the ABC are, the better a deviation from this specific angle should be taken into account. 
Therefore, it is worth studying how the ABCs
react to different angles $\theta$ of incidence. In this respect, we consider $100\, {\rm kHz}$ waves traveling from left to right in a 
rectangular waveguide $\Omega=[0,3\, {\rm cm}]\times[0, 0.05\, {\rm cm}]$ on the right wall of which one of the studied ABCs is set.
We begin from showing the well-known effect: the more the deviation of the 
incidence angle from zero the more the solution is contaminated by reflected waves.
The first example presented in this series is a wave impinging the boundary at $\theta=0^{\circ}$. In this simple case all the ABCs should work the best.
The results are shown in figure~\ref{fig:ABC_2d_1}.
As in the 1-d case, $\rm ABC^{1,1}_{PS}$ outperforms $\rm ABC^{1,0}_{PS}$. However,
$\rm ABC^{2,2}_{EM}$, which coincides with $\rm ABC^{1,1}_{EM}$ for $\theta=0^{\circ}$, shows low accuracy.

\begin{figure}[htp]
  \hspace*{1.3cm}
  \begin{tabular}{cc}
  \hspace*{0.5cm}\large$\delta$ &
  \begin{minipage}{0.9\textwidth}
    \includegraphics[width=3in]{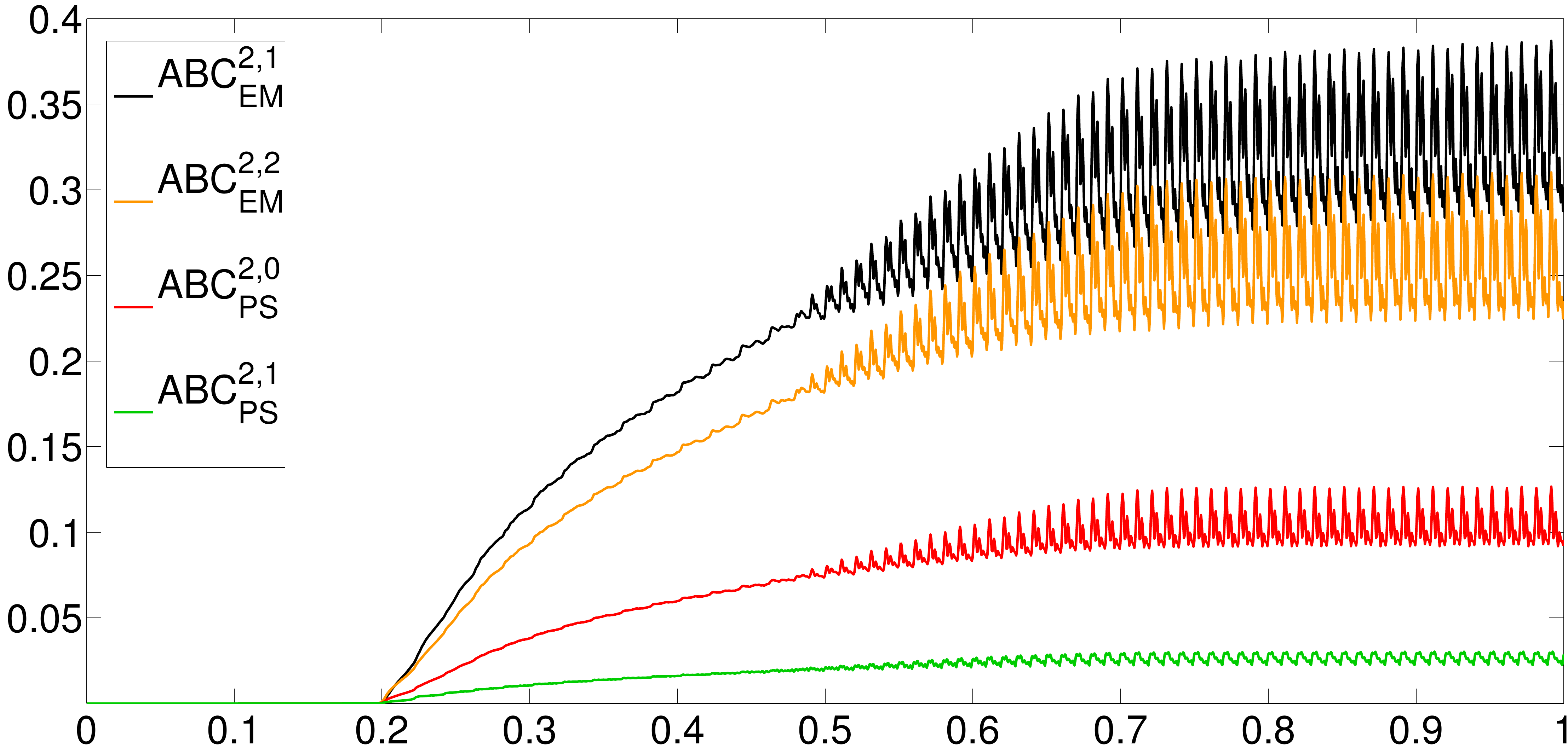}\\
  \end{minipage}\\ 
  \vspace*{-0.7cm}&\\
  & \hspace{-3.4cm}\large $t$  \\
  \end{tabular}  
  \caption{Relative error $\delta$ versus time $t$ for $\rm ABC^{2,o}_n$, $\rm o=\{0,1,2\}$, $\rm n=\{EM,PS\}$.}
  \label{fig:ABC_2d_2}
\end{figure}

In the next setup, we increase the angle $\theta$ by $15^{\circ}$, which is a assumed to be a harder trial for the ABCs (see figure ~\ref{fig:ABC_2d_2}).
Indeed, all the ABCs revealed to be quite sensitive to the angle of incidence and exhibit higher errors introduced
by the reflected waves into the solution. The second order Engquist--Majda ABC gives more accurate results
compared to the first order condition, but the error is quite large.
The new boundary condition of the first order demonstrates the lowest error while the zero order ABC is less efficient.

\begin{figure}[htp]
  \hspace*{1.3cm}
  \begin{tabular}{cc}
  \hspace*{0.5cm}\large$\delta$ &
  \begin{minipage}{0.9\textwidth}
    \includegraphics[width=3in]{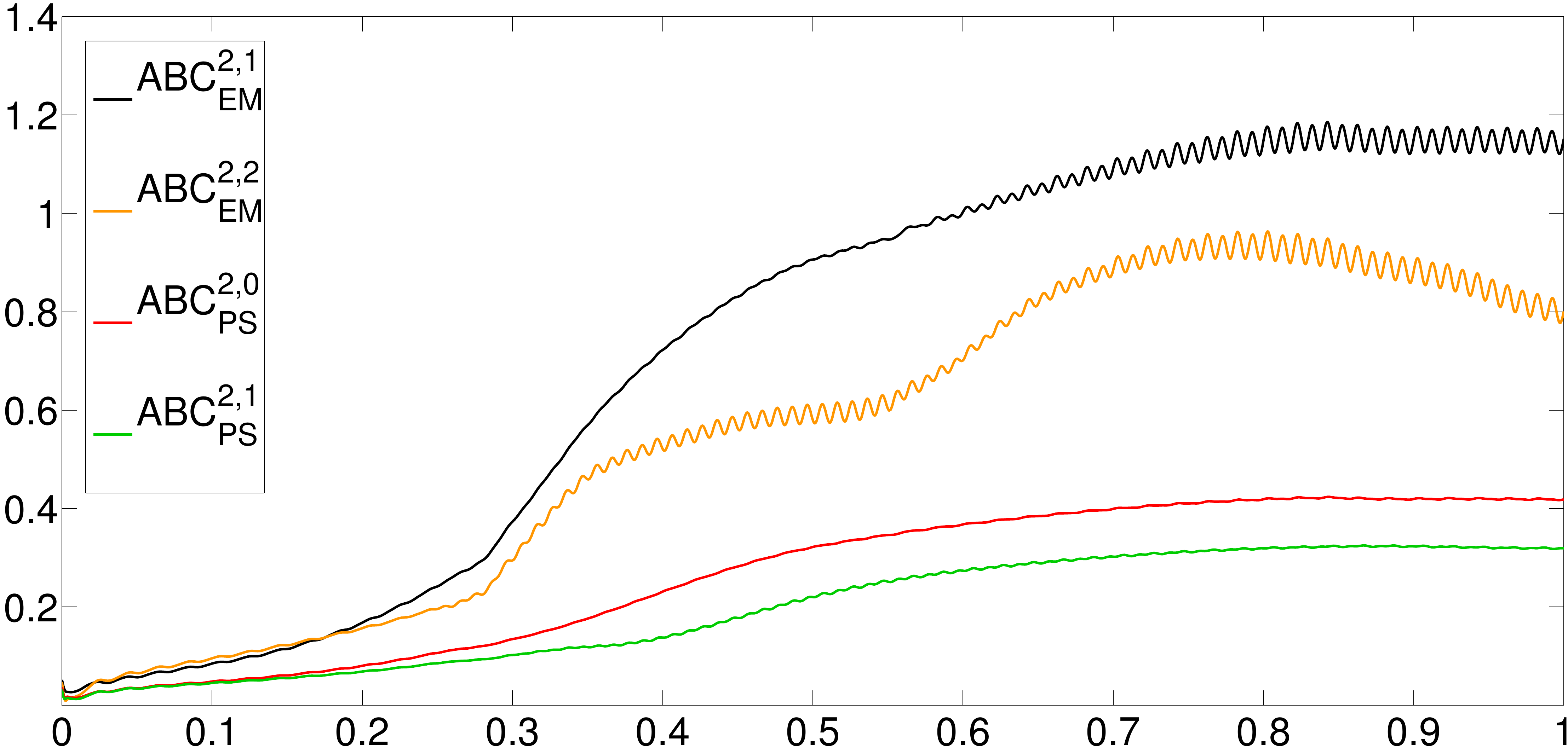}\\
  \end{minipage}\\ 
  \vspace*{-0.7cm}&\\
  & \hspace{-3.5cm}\large $t$  \\
  \end{tabular}  
  \caption{Relative error $\delta$ versus time $t$ for $\rm ABC^{2,o}_n$, $\rm o=\{0,1,2\}$, $\rm n=\{EM,PS\}$.}
  \label{fig:ABC_2d_3}
\end{figure}

In the last example, the ABCs are studied in a much more realistic situation -- the HIFU problem which is routinely used in 
computational setups to simulate the thermotherapy for human liver cancer. 
We consider a concave transducer, with a much higher excitation frequency $\omega=1\, {\rm MHz}$, located
at bottom of a square domain $\Omega=[0,20\, {\rm mm}]\times[0, 20\, {\rm mm}]$. The numerical results are given in figure~\ref{fig:ABC_2d_3}.

At the very beginning of the simulation ($t<0.3$), the
first and second order Engquist--Majda ABCs work equally well. However, the situation gets worse as time advances:
the discrepancy between the boundary conditions grows and the solution becomes substantially contaminated by the reflected waves. 
The second order Engquist--Majda ABC does not dramatically affect the situation, and the numerical
solution is still quite poor. The proposed ABCs demonstrate much lower errors. However, the difference
between them is less pronounced compared to the waveguide example.
Another remarkable feature of the new ABC is that the error exhibit a much less fluctuating behavior, which suggests that
the new boundary conditions are robust with respect to the wave propagation regime. 

\section{Conclusions}
In this work we proposed zero and first order ABCs, based on pseudo-differential calculus,
for the Westervelt equation in one and two space dimensions. 
Well-posedness of the boundary value problem with the new ABCs is stated and proven.
All our numerical results reconfirm the fact that using the ABCs which are not especially tailored for the Westervelt equation
lead to poor numerical solutions. The zero order ABCs are computationally easier than the first order conditions, however,
more prone to the regimes of the wave propagation and less accurate. It is important to remark that the application
of the self-adapting technique~\cite{SW2012} to the developed ABCs will result in further improvements.

\section*{Acknowledgments}
The first author would like to thank EPSRC Mathematics Platform grant EP/I019111/1, which partly
supported this work.
The second author gratefully acknowledges support by the Austrian Science Fund (FWF) under the grant P24970. Moreover we thank Barbara Wohlmuth, TU Munich, and Manfred Kaltenbacher, TU Vienna, for stimulating discussions

\subsection{Appendix: ABCs for the 1-d Westervelt equation via a paradifferential approach }
\label{sec_paradiff}
In this section we focus on the construction of transparent boundary conditions directly for the nonlinear Westervelt equation.
The disadvantage of the pseudo-differential approach for designing ABCs is in its inability to treat nonlinear equations. This
obstacle can be overcome by using the para-differential calculus originated from the paper of Bony~\cite{Bony1981} 
with an improvement done by Meyer~\cite{Meyer1981}. Despite the para-differential calculus and especially the para-linearization 
technique of Bony embrace wide opportunities to build ABCs for nonlinear equations, their use is very restricted in current research works.
The first application of para-differential operators to the development of ABCs was done for the Burgers equation in~\cite{Dubach2000}.
Some relatively recent results can be found in few works (see~\cite{Szeftel2006_1,Szeftel2006_2}).

Before the derivation of ABCs we briefly recall some general facts about para-differential operators and Bony's para-linearization.
In according to~\cite{Bony1981}, the multiplication by a function $a(x)\in C^{\infty}$, $x\in\mathbb{R}^d$ can
be replaced with the operator $T_a$ defined as
\begin{equation}
 \mathcal{F}(T_au)(\zeta)=\frac{1}{(2\pi)^d}\int\limits_{\mathbb{R}^d}\chi(\zeta-\eta,\eta)\mathcal{F}a(\zeta-\eta)\mathcal{F}u(\eta)\ d\eta,
 \label{BonyParaProduct}
\end{equation}
where $\mathcal{F}$ is the Fourier transformation, $\chi\in C^{\infty}(\mathbb{R}^d\times\{\mathbb{R}^d\setminus\{0\}\})$ is a function
of homogeneity degree zero satisfying
\begin{equation}
\left\{
\begin{aligned}
\chi(\zeta,\eta)=1&\quad \text{if }\ |\zeta|\le\varepsilon_1|\eta|,\\
\chi(\zeta,\eta)=0&\quad \text{if }\ |\zeta|\ge\varepsilon_2|\eta|,
\end{aligned}
 \right.
\end{equation}
where $0<\varepsilon_1<\varepsilon_2$.

Let us consider a nonlinear differential equation of order $N$ defined by the superposition operator induced by $\Phi$
\begin{equation}
F[u](x)=\Phi(x,u(x),\ldots,\partial^{\alpha}u(x),\ldots)_{0\le|\alpha|\le N}=0 
\label{NonlinEqF}
\end{equation}
with $\Phi\in C^{\infty}$ and $x\in\mathbb{R}^d$.
In accordance to~\cite{Bony1981}, the para-linearization of~(\ref{NonlinEqF}) 
with $\Phi(x,\cdot)$ vanishing at $0$ is given by
\begin{equation}
 F[u]=\sum_{0\le|\alpha|\le N} T_{\frac{\partial \Phi}{\partial \lambda_\alpha} (\cdot,u,\ldots,\partial^{\alpha}u,\ldots)_{0\le|\alpha|\le N}} \partial^{\alpha}u +R(u),
 \label{Paralin_F}
\end{equation}
where $T_{a}$ is a para-differential operator with symbol $a$, and $R(u)$ is a smooth error. More precisely,
for all $u\in H^s(\mathbb{R}^d)$ with $s>d/2$ equation~(\ref{Paralin_F}) implies $R(u)\in H^{2s-d/2}$ (see~\cite{Meyer1981}). Equation~(\ref{Paralin_F}) is often referred to as
the para-linearization formula of Bony.

Before the derivation of ABCs for the Westervelt equation~(\ref{Westervelt_1D}), we develop
the nonlinear term on its right hand side $\gamma(u^2)_{tt}=2\gamma((u_t)^2+uu_{tt})$
and recast~(\ref{Westervelt_1D}) in the form
\begin{equation}
\nu^2(u)u_{tt}-u_{xx}-\beta u_{txx}=2\gamma(u_t)^2
\label{Westervelt_1D_2_}
\end{equation}
with $\nu^2(u)=c^{-2}-2\gamma u$.

Based on~(\ref{Paralin_F}) and taking into account that (see~\cite{Bony1981})
\begin{equation}
  fg=T_fg+T_gf+R,
  \label{BonyParaProduct1}
\end{equation}
where $T_f$ and $T_g$ are para-differential operators with symbols $f$ and $g$, we get a para-differential equation

\begin{equation}
\mathfrak{D}_2u=0,\quad \mathfrak{D}_2=c^{-2}\partial^2_t-2\gamma(T_{u_{tt}}+T_u\partial^2_t)-\partial^2_x-\beta\partial_{txx}-2\gamma T_{2u_t}\partial_t
\label{Westervelt_1D_2__}
\end{equation}
instead of the nonlinear Westervelt equation~(\ref{Westervelt_1D}).

Acting similar to the previous derivation, we can apply Nirenberg's factorization, 
analogous to~(\ref{NFact_Westervelt_1D_u1}), and rewrite~(\ref{Westervelt_1D_2__}) in the form
\begin{equation}
\mathfrak{D}_2=-(\partial_x-A)(\partial_x-B)+R,
\label{NFact_Westervelt_1D_u2}
\end{equation}
where $A$ and $B$ are para-differential operators with symbols $a$ and $b$, respectively.

A similar argument as for the linearized Westervelt equation 
yields
\begin{equation}
\nu^2(u)(i\tau)^2-2\gamma u_{tt}-\beta(i\tau)\partial^2_x-4\gamma u_t(i\tau)=(a+b)\partial_x+\partial_xb-ab+R.
\label{NFact_Westervelt_1D_sym1}
\end{equation}
Note that this equation differs from \eqref{NirFact_D1_sym1} and it will also lead to different ABCs. 
Again we will skip the $\beta$ terms for the same reason as in section \ref{sec:LinWestervelt_1d}

The equation for the $\mathcal{O}(\tau^2)$ terms remains the same, namely \eqref{a1b1_eq} 
Thus we get 
\begin{equation}
b_1=-a_1=\nu(i\tau)
\label{a1b1_1D_1_para} 
\end{equation}

The equation for the $\mathcal{O}(\tau^1)$ terms is different from \eqref{a0b0_1D_tau1}, namely
\begin{equation}
 \left\{
  \begin{aligned}
   a_0+b_0&=0,\\
   \beta(i\tau)\partial^2_x+4\gamma u^{(0)}_t(i\tau)&=a_1b_0+a_0b_1-i{a_1}_{\tau}{b_1}_t-{b_1}_x,
  \end{aligned}
  \right.
\label{a0b0_1D_tau1_para}
\end{equation}
which upon setting $\beta=0$ yields (differently from \eqref{a0b0_1D_eq2_5})
\begin{equation}
b_0=-a_0=-\frac{1}{2\nu}\left(\mathcal{A}_0[\nu]
+4\gamma u_t\right)
\label{a0b0_1D_para}
\end{equation}
with the operator $\mathcal{A}_0:=\partial_x+\nu\partial_t$. 

Finally the $\mathcal{O}(\tau^0)$ equation is (differently from \eqref{am1bm1_1D_tau0_1})
\begin{equation}
 \left\{
  \begin{aligned}
   a_{-1}+b_{-1}&=&0,\\
   -a_1b_{-1}-a_0b_0-a_{-1}b_1+i({a_1}_{\tau}{b_0}_t+{a_0}_{\tau}{b_1}_t)-\frac{i^2}{2}{a_1}_{\tau\tau}{b_1}_{tt}+{b_0}_x&=&-2\gamma u_{tt},
  \end{aligned}
  \right.
\label{am1bm1_1D_tau0_para}
\end{equation}
so that we get
\begin{equation}
\begin{split}
b_{-1}=-a_{-1}=\frac{1}{2\nu(i\tau)}&
\left(\mathcal{A}_0\left[\frac{1}{2\nu}\left(\mathcal{A}_0[\nu]
+4\gamma u_t\right)\right]\right.
\left. -\left(\frac{1}{2\nu}\left(\mathcal{A}_0[\nu]
+4\gamma u_t\right)\right)^2
-2\gamma u_{tt}\right) =:\frac{\tilde{\mu}}{2\nu(i\tau)}.
\end{split}
\label{am1bm1_1D_tau0_3}
\end{equation}

Parallel to the ABCs for the linearized Westervelt equation from section \ref{sec:LinWestervelt_1d}, we obtain the zero, first and second order boundary conditions
\begin{equation}
\left.\mathcal{A}^{\prime}_0u\right|_{x=0}=\left.\left(\partial_x+\nu(u)\partial_t\right)\right|_{x=x_{\rm ABC}} u=0
\label{ParaDiffABC0_1D_xa}
\end{equation}
\begin{equation}
\left.\mathcal{A}^{\prime}_1u\right|_{x=0}=\left.(\mathcal{A}^{\prime}_0-\mathcal{B}^{\prime}_1)u\right|_{x=0}=
\left.\left(\partial_x+\nu(u)\partial_t-\frac{1}{2\nu(u)}\left(\mathcal{A}^{\prime}_0[\nu(u)]+4\gamma u_t\right)\right)u \right|_{x=0} =0
\label{ParaDiffABC1_1D_xa}
\end{equation}
with $\mathcal{B}^{\prime}_1:=\frac{1}{2\nu(u)}(\mathcal{A}^{\prime}_0[\nu]+4\gamma u_t))$,
\begin{equation}
\left.\mathcal{A}^{\prime}_2u\right|_{x=0}=\left.\left(\mathcal{A}^{\prime}_1 u_t -\mathcal{B}^{\prime}_2 u\right)\right|_{x=0}=
\left.\left(u_{xt}+\nu u_{tt}-\frac{1}{2\nu}\left((\nu_x+\nu \nu_t)u_t+4\gamma (u_t)^2-\tilde{\mu} u\right)\right)
\right|_{x=0}=0\,,
\label{ParaDiffABC2_1D_xa}
\end{equation}
where $\mathcal{B}^{\prime}_2:=\tilde{\mu}(u)$, which contains multiplication with $u_{tt}$, as opposed to \eqref{ABC2_1D_xa}.

As in the pseudo-differential case, we do not consider higher order boundary conditions, although their derivation follows the same lines.

\begin{remark}
It is probably due to the approximation by Taylor linearization \eqref{Paralin_F} for an actually quadratic nonlinearity that the paradifferential approach yields a different and in our view most likely worse ABC than the pseudodifferential one from section \ref{sec:ABCWestervelt}, since the linearization \eqref{Westervelt_lin} seems to be better capable to capture the quadratic nonlinearity than just Taylor linearization (cf. \eqref{Westervelt_linTaylor}).
\end{remark}

\end{document}